\documentclass[a4paper,preprint,10pt,3p,numbers,sort&compress]{elsarticle}

\usepackage{cancel}
\usepackage{comment}
\usepackage[ruled,vlined]{algorithm2e}
\usepackage[english]{babel}
\usepackage{hyperref}
\usepackage{float}
\usepackage{graphicx}
\usepackage{amssymb}
\usepackage{amsmath}
\usepackage{tikz}
\usepackage{standalone}
\usepackage{booktabs}

\usepackage[utf8]{inputenc}
\usepackage{pgfplots}
\DeclareUnicodeCharacter{2212}{−}
\usepgfplotslibrary{groupplots,dateplot}
\usetikzlibrary{patterns,shapes.arrows}
\pgfplotsset{compat=newest}

\usepackage{caption}
\usepackage{subcaption}

\begin{document}

\begin{frontmatter}

\title{Energy-Conserving Neural Network for Turbulence Closure Modeling} 
\author[1]{T. van Gastelen}
\author[1]{W. Edeling}
\author[1]{B. Sanderse}


\address[1]{Centrum Wiskunde \& Informatica, 
            Science Park 123, 
            Amsterdam,
            The Netherlands}


\begin{abstract}
In turbulence modeling, we are concerned with finding closure models that represent the effect of the subgrid scales on the resolved scales. Recent approaches gravitate towards machine learning techniques to construct such models. However, the stability of machine-learned closure models and their abidance by physical structure (e.g. symmetries, conservation laws)  are still open problems. To tackle both issues, we take the `discretize first, filter next' approach. In this approach we apply a spatial averaging filter to existing fine-grid discretizations. The main novelty is that we introduce an additional set of equations which dynamically model the energy of the subgrid scales. Having an estimate of the energy of the subgrid scales, we can use the concept of energy conservation to derive stability. The subgrid energy containing variables are determined via a data-driven technique. The closure model is used to model the interaction between the filtered quantities and the subgrid energy. Therefore the total energy should be conserved. Abiding by this conservation law yields guaranteed stability of the system. In this work, we propose a novel skew-symmetric convolutional neural network architecture that satisfies this law. The result is that stability is guaranteed, independent of the weights and biases of the network. Importantly, as our framework allows for energy exchange between resolved and subgrid scales it can model backscatter. To model dissipative systems (e.g. viscous flows), the framework is extended with a diffusive component. The introduced neural network architecture is constructed such that it also satisfies momentum conservation. We apply the new methodology to both the viscous Burgers' equation and the Korteweg-De Vries equation in 1D. The novel architecture displays superior stability properties when compared to a vanilla convolutional neural network.
    
\end{abstract}

\end{frontmatter}

\section{Introduction}

Direct numerical simulations (DNSs) of turbulent flows are often infeasible due to the high computational requirements. Especially for applications in design and uncertainty quantification this rapidly becomes computationally infeasible, as typically many simulations are required \cite{smith2013uncertainty,sasaki2002navier}. To tackle this issue several different approaches have been proposed, such as reduced order models \cite{podbenjamin}, Reynolds-averaged Navier-Stokes (RANS) \cite{RANS}, and Large Eddy Simulation (LES) \cite{sagaut2006large}. These approaches differ in how much of the physics is modeled. Here we will focus on the LES approach. 

In LES, the large-scale physics is modeled directly by a coarse-grid discretization. The coarse grid is accounted for by applying a filter to the original equations. However, as the filter does not commute with the nonlinear terms in the equations, a commutator error arises. This prevents one from obtaining an accurate solution without knowledge of the subgrid-scale (SGS) content. This commutator error is referred to as the closure term. Modeling this term is the main concern of the LES community. A major difficulty in this process is dealing with energy moving from the small scales to the large scales (backscatter) \cite{backscatter1,backscatter2}. This is because the SGS energy is unknown during the time of the simulation. This makes accounting for backscatter without leading to numerical instabilities difficult \cite{instability1}. Classical physics-based closure models are therefore often represented by a dissipative model, e.g.\ of eddy-viscosity type \cite{smagorinsky1963general}. This ensures a net decrease in energy. Another option is that the closure model is clipped such that backscatter is removed \cite{clipping}. Even though the assumption of a global net decrease in energy is valid \cite{smagorinsky1963general}, explicit modeling of backscatter is still important. This is because locally the effect of backscatter can be of great significance \cite{backscatter_importance1,backscatter_importance2}. Closure models that explicitly model the global SGS energy at a given point in time, to allow for backscatter without sacrificing stability, also exist \cite{ghosal_total_energy}. Recently, machine learning approaches, or more specifically neural networks (NNs), have come forward as viable closure models. They have been shown to outperform the classical approaches for different use cases \cite{Frezat,List,park,guan,kochkov}. However, stability remains an important issue, along with abidance by physical structure such as mass, momentum, and energy conservation \cite{park,beck1,beck2,beck3}.

In \cite{beck1} homogeneous isotropic turbulence for the compressible Navier-Stokes equations was treated. A convolutional neural network (CNN) was trained to reproduce the closure term from high-resolution flow data. Although \textit{a priori} cross-correlation analysis on the training data showed promising results, stable models could only be achieved by projecting onto an eddy-viscosity basis. In \cite{beck2} a gated recurrent NN was applied to the same test case. This network displayed even higher cross-correlation with the closure term, but still yielded unstable models. Even after training on data with added artificial noise the model remained unstable \cite{beck3}. In \cite{park} incompressible turbulent channel flow was treated. Here NNs with varying levels of locality were used to construct a closure model. They showed that increasing the view of the NN improves \textit{a priori} performance. However, \textit{a posteriori} analysis showed that this increased input space also led to instabilities. Even after introducing backscatter clipping these larger models were still outperformed by the highly localized NN models. Two promising approaches to improving the stability of NN closure models are `trajectory fitting' \cite{Frezat,List,MacArt_embedded_learning,Syver,Melchers2022} and reinforcement learning \cite{beck4,Bae2022}. Both of these approaches have in common that instead of fitting the NN to the exact closure term (which is what we will refer to as `derivative fitting'), one optimizes directly with respect to how well the solution is reproduced. This has been shown to yield more accurate and stable closure models \cite{Frezat,List,Melchers2022}. The difference between these two methods is that trajectory fitting uses exact gradients, such that gradient-based optimizers can be applied to optimize the NN weights \cite{kingma_adam}. Reinforcement learning does not require these gradients. This makes it suitable for non-differentiable processes such as chess and self-driving cars \cite{reinforcement_learning}.

These approaches have all been applied to the Navier-Stokes equations. However, in this paper we consider a 1D simplification, namely Burgers' equation. Several studies have been carried out which apply machine learning to this equation. In \cite{burgers_pinn} physics-informed neural networks (PINNs) were succesfully applied to Burgers' equation. In PINNs the PDE is encoded into the loss function of the neural network. The advantage of PINNs is that they allow us to approximate the solution without explicitly discretizing space and time. In the context of closure modeling, several studies have been carried out. For example, in \cite{Bar-Sinai} and \cite{burgers_2} a neural network was fitted to predict the closure term for Burgers' equation with forcing. In \cite{burgers_2} the trained neural network was successfully applied to unseen viscosity values. This was achieved by limited retraining on new data. This technique is known as transfer learning. Furthermore, several studies have been carried out which show the benefits of trajectory fitting \cite{Melchers2022,burgers_3,burgers_1}. In \cite{burgers_3} trajectory fitting is combined with a Fourier neural operator to predict a spatially dependent Smagorinsky coefficient. This combination outperformed the other considered approaches. It also guarantees stability by being strictly dissipative. Furthermore, Fourier neural operators have the advantage that they are grid independent \cite{fourier_neural_operators}.

However, none of the discussed approaches leads to a provably stable NN closure model, while still allowing for backscatter. In addition, they do not guarantee abidance by the underlying energy conservation law. The latter is something that to our knowledge does not yet exist in the case of LES closure models. To resolve these shortcomings, we present \textit{a new NN closure model that satisfies both momentum and kinetic energy conservation and is therefore stable by design}. As stated earlier, the difficulty of this task mainly lies in the fact that: (i) The kinetic energy conservation law includes terms which depend on the SGS content which is too expensive to simulate directly. (ii) Consequently, the kinetic energy of the large scales is not a conserved quantity (in the limit of vanishing viscosity). To tackle these issues we take the `discretize first, filter next' approach \cite{Syver,Melchers2022}. This means that we start from a high-resolution discretization with $N$ degrees of freedom, to which we apply a discrete filter. This filter projects the solution onto a coarse computational grid of dimension $I$, with $I\ll N$. Given the discrete filter, the exact closure term can be obtained by computing the commutator error. The main advantage of this approach is that the closure term now also accounts for the discretization error. Based on the filter's properties we then derive an energy conservation law that can be split into two components: one that depends solely on the resolved scales (resolved energy) and another that solely depends on the SGS content (SGS energy) \cite{ghosal_total_energy}. Like in existing works the closure model is represented by a convolutional neural network (CNN) \cite{conv_NN}. The main novelty comes from the addition of a set of SGS variables. These SGS variables represent the SGS energy, projected onto the coarse grid.
The key insight is that the resulting system of equations should still conserve energy in the inviscid limit. We then choose our CNN architecture such that it is consistent with this limit. In this way we still allow for backscatter without sacrificing stability.

The paper is structured in the following way: In section \ref{sec:preliminaries} we discuss Burgers' and Korteweg-de Vries equation and their energy and momentum conservation properties. We introduce the discrete filter, the resulting closure problem, and derive a new energy conservation law. This law describes the exchange between the resolved and the SGS energy. In section \ref{sec:framework} we introduce our novel machine learning framework for modeling the closure term. This approach satisfies the derived energy conservation law using the set of SGS variables to represent the SGS energy. In addition, we show how to satisfy momentum conservation. In section \ref{sec:results} we study the convergence and stability of our closure modeling framework and compare this to a vanilla CNN. We also analyze its structure-preserving properties in terms of momentum and energy conservation and its ability to extrapolate in space and time. In section \ref{sec:discussionNS} we present a short discussion on the applicability of our framework to the Navier-Stokes equations. In section \ref{sec:conclusion} we conclude our work.

\section{Governing equations, discrete filtering, and closure problem}\label{sec:preliminaries}

Before constructing a machine learning closure, we formulate a description of the closure problem on the discrete level. For this purpose we introduce the filter and reconstruction operator which we apply to the discrete solution. In this way we also account for the discretization error. In addition, we discuss the effects of filtering on the physical structure of the system. 

\subsection{Spatial discretization}
We consider an initial value problem of the following form:
\begin{align}
    \frac{\partial u}{ \partial t} &= f(u), \label{eq:IVP_continuous} \\
    u(\mathbf{x},0) &= u_0(\mathbf{x}),
\end{align}
which describes the evolution of some quantity $u(\mathbf{x},t)$ in space $\mathbf{x} \in \Omega$ and time $t$ on the spatial domain $\Omega \subseteq \mathbb{R}^d$, given initial state $u_0$. The dynamics of the system is governed by the right-hand side (RHS) $f(u)$, which typically involves partial derivatives of $u$ with respect to $\mathbf{x}$. After spatial discretization (method of lines), we obtain the vector $\mathbf{u}(t) \in \mathbb{R}^N$. The elements $\text{u}_i$ of this vector approximate the value of $u$ at each of the $N$ grid points $\mathbf{x}_i \in \Omega$ for $i = 1,\ldots,N$, such that $\text{u}_{i} \approx u(\mathbf{x}_i)$. The discrete analogue of the IVP is then
\begin{align}\label{eq:discr_ODE}
     \frac{\text{d} \mathbf{u}}{ \text{d} t} &= f_h(\mathbf{u}), \\
    \mathbf{u}(0) &= \mathbf{u}_0,
\end{align}
where $f_h$ represents a spatial discretization of $f$. It is assumed that all the physics described by equation \eqref{eq:IVP_continuous} is captured in the discrete solution $\mathbf{u}$. This means that whenever the physics involves a wide range of spatial scales, a very large number of degrees of freedom $N$ is needed to adequately resolve all the scales. This results in a large amount of computational resources required to solve these equations. This is what we aim to alleviate.

\subsection{Burgers' and Korteweg-de Vries equation, and physical structure}\label{sec:burgers}

We are interested in the modeling and simulation of turbulent flows. For this purpose, we first consider Burgers' equation, a 1D simplification of the Navier-Stokes equations.
Burgers' equation describes the evolution of the velocity $u(x,t)$ according to partial differential equation (PDE)
\begin{align}\label{eq:burgerseq}
    \frac{\partial u}{\partial t} &= -\frac{1}{2}\frac{\partial u^2}{\partial x} + \frac{\partial}{\partial x} \left( \nu\frac{\partial u}{\partial x} \right).
\end{align}
The first term on the RHS represents nonlinear convection and the second term diffusion, weighted by the positive viscosity field $\nu(x) \geq 0$. 
This equation expresses similar behavior to 3D turbulence in the fact that smaller scales are created by the nonlinear convective term, which then dissipate through diffusion \cite{love_1980}. We will be interested in two properties of the Burgers' equation, which we collectively call `structure'. 

Firstly, momentum $P$ is conserved on periodic domains:
\begin{equation}\label{eq:continuous_mom_cons}
\frac{\text{d}P}{\text{d}t} =  \frac{\text{d}}{\text{d}t} \underbrace{\int_{\Omega}  u \text{d}\Omega}_{=:P} =   \int_{\Omega} -\frac{1}{2}\frac{\partial u^2}{\partial x} + \frac{\partial}{\partial x} \left( \nu\frac{\partial u}{\partial x} \right) \text{d}\Omega = 0.   
\end{equation}
Secondly, on periodic domains (kinetic) energy $E$ is conserved in the absence of viscosity: 
\begin{equation}\label{eq:continuous_E_cons}
  \frac{\text{d} E}{\text{d} t} = \frac{\text{d}}{\text{d} t}\underbrace{\frac{1}{2}\int_{\Omega}u^2\text{d}\Omega}_{=:E} 
  = \int_{\Omega}-  \frac{u}{2} \frac{\partial u^2}{\partial x} + u \frac{\partial}{\partial x} \left( \nu\frac{\partial u}{\partial x} \right)\text{d}\Omega  = -  \int_\Omega \nu \left(\frac{\partial u}{\partial x}\right)^2\text{d}\Omega \leq 0,
\end{equation}
where we used integration by parts. As we solely deal with viscous flows, i.e.\ $\nu(x) > 0$, we disregard energy loss due to the presence of shocks \cite{Jameson}. Note that these conservation laws only hold in the absence of forcing.

These properties can be preserved in a discrete setting by employing a structure-preserving scheme \cite{Jameson}. On a uniform grid the convective term is discretized with the following skew-symmetric scheme:
\begin{equation}
   (\mathbf{C}(\mathbf{u})\mathbf{u})_i =  -\frac{1}{3h}(\text{u}^2_{i+1}-\text{u}^2_{i-1}) - \frac{1}{3h}\text{u}_{i} (\text{u}_{i+1} - \text{u}_{i-1}),
\end{equation}
where $h$ is the grid spacing. The skew-symmetry entails that $\mathbf{u}^T\mathbf{C}(\mathbf{u})\mathbf{u}=0$. This is used later to derive energy conservation.
Furthermore, the diffusion operator is discretized as 
\begin{equation}\label{eq:Q^TQ}
    (-\mathbf{Q}^T\text{diag}(\boldsymbol{\nu})\mathbf{Q} \mathbf{u})_{i} = \frac{1}{h^2}(\nu_{i}(\text{u}_{i+1} - \text{u}_{i})+ \nu_{i-1}(\text{u}_{i-1} -\text{u}_{i})),
\end{equation}
where $\mathbf{Q}$ is a simple forward difference approximation of the first derivative and $\nu_i = \nu(x_i)$ \cite{podbenjamin,benjamin_thesis}.
In this paper we deal with constant viscosity $\nu(x)=\nu$ such that we obtain the following system of ordinary differential equations (ODEs):
\begin{equation}\label{eq:discr_burgers}
    \frac{\text{d}\mathbf{u}}{\text{d}t} = \mathbf{C}(\mathbf{u})\mathbf{u} + \nu \mathbf{D}\mathbf{u}.
\end{equation}
Here $\mathbf{D}=-\mathbf{Q}^T\mathbf{Q}$ corresponds to a simple central difference approximation of the second derivative.
For the time discretization we employ an explicit RK4 scheme \cite{RK4_Butcher:2007}.

This discretization conserves the discrete momentum $P_h = h\mathbf{1}^T\mathbf{u}$ in the periodic case:
\begin{equation}
    \frac{\text{d}P_h}{\text{d}t} = h\mathbf{1}^T \frac{\text{d}\mathbf{u}}{\text{d}t} = 0,
\end{equation}
where $\mathbf{1}$ is a column vector with all entries equal to one. Note that this equation discretely represents the integral in \eqref{eq:continuous_mom_cons}. Furthermore, due to the skew-symmetry of the convection operator the evolution of the discrete kinetic energy $E_h = \frac{h}{2}\mathbf{u}^T\mathbf{u}$ is given by:
\begin{equation}\label{eq:DNS_stability_cond_burgers}
    \text{Burgers' equation:}\qquad \frac{\text{d}E_h}{\text{d}t} = h\mathbf{u}^T \frac{\text{d}\mathbf{u}}{\text{d}t}=h\nu\mathbf{u}^T\mathbf{D}\mathbf{u} = - h \nu  ||\mathbf{Q}\mathbf{u}||_2^2 \leq 0.
\end{equation}
This is the discrete equivalent of \eqref{eq:continuous_E_cons}. In both the continuous and discrete formulation we used the product rule to obtain the derivative. The norm $\| .\|_{2}$ represents the conventional 2-norm. From \eqref{eq:DNS_stability_cond_burgers} we conclude that this discretization ensures net kinetic energy dissipation, and conservation in the inviscid limit. From this point forward we will refer to the kinetic energy simply as energy.

In addition to Burgers' equation we will consider the Korteweg-de Vries (KdV) equation:
\begin{equation}\label{eq:KdV}
    \frac{\partial u}{\partial t} = - \frac{\varepsilon}{2}\frac{\partial u^2}{\partial x} - \mu\frac{\partial^3 u}{\partial x^3},
\end{equation}
where $\varepsilon$ and $\mu$ are scalar parameters.
The KdV equation conserves momentum and energy irrespective of the values of $\varepsilon$ and $\mu$. We discretize the nonlinear term in the same way as for Burgers' equation, using the skew-symmetric scheme. The third-order spatial derivative is approximated by a skew-symmetric central difference stencil: $(-\text{u}_{i-2}+2\text{u}_{i-1}-2\text{u}_{i+1}+\text{u}_{i+2}) /(2h^3)$, see  \cite{KdV_Yan2019}. The resulting discretization is not only momentum conserving, but also energy conserving:
\begin{equation}\label{eq:DNS_stability_cond_KdV}
    \text{KdV equation:}\qquad  \frac{\text{d}E_h}{\text{d}t} = 0.
\end{equation}

\subsection{Discrete filtering}\label{sec:discrete_filtering}

In order to alleviate the high computational expenses for large $N$ we apply a spatial averaging filter to the fine-grid solution $\mathbf{u}$. This results in the coarse-grid approximation $\bar{\mathbf{u}} \in \mathbb{R}^I$.  The coarse grid follows from dividing  $\Omega$ into $I$ non-overlapping cells $\Omega_i$ with cell centers $\mathbf{X}_i$. The coarse grid is refined into the fine grid by splitting each $\Omega_i$ into $J$ subcells $\omega_{ij}$ with cell centers $\mathbf{x}_{ij}$.
The subdivision is intuitively pictured in Figure \ref{fig:grid}, for a uniform 1D grid. Furthermore, we define the mass matrices $\boldsymbol{\omega} \in \mathbb{R}^{N\times N}$ and $\boldsymbol{\Omega} \in \mathbb{R}^{I \times I}$ which contain the volumes of the fine and coarse cells on the main diagonal, respectively.

\begin{figure}[ht]
    \centering
    \begin{tikzpicture}
    \draw  (1,-0.25) -- (1,0.25);
    \draw  (1,-1) -- (1,0.25-1);
    \foreach \i in {1,2,3}{
    \draw  (3*\i+1,-0.25) -- (3*\i+1,0.25);
    
    \draw  (3*\i+1,-1) -- (3*\i+1,0.25-1);
    
    \draw  (1+3*\i-3,-1) -- (1+3*\i,-1);
    \fill[black]  (1+3*\i-1.5,-1) circle (0.07cm) node[above]{$ \bar{\text{u}}_{\i}$} node[below]{$ \Omega_{\i }$} ;

        \foreach \j in {1,2,3}{
        \def\ij{\the\numexpr 3 * \i - 3 + \j}
        \draw  (1+3*\i-3+\j,-0.125) -- (1+3*\i-3+\j,0.125);
        \fill[black]  (1+3*\i-3+\j -0.5,0) circle (0.07cm) node[above]{$\mathsmaller{ \text{u}_{\ij}}$} node[below]{$\mathsmaller{ \omega_{\i \j}}$};
        \draw (3*\i-3+\j,0) -- (1+3*\i-3+\j,0);
    
        }}
        
    \end{tikzpicture}
    \caption{Subdivision of the spatial grid. The dots represent cell centers $\mathbf{x}_{ij}$ and $\mathbf{X}_{i}$ for $N=9$ and $I=J=3$.}
    \label{fig:grid}
\end{figure}

\begin{figure}[ht]
    \centering
    \includegraphics[width = 0.48\textwidth]{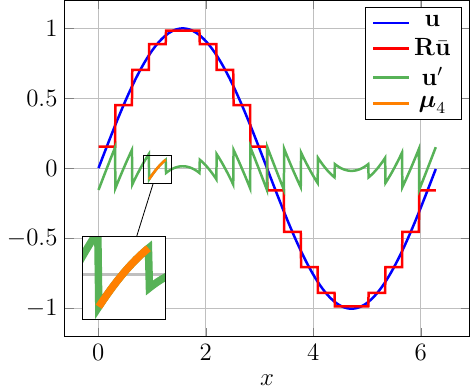} \hfill
    \includegraphics[width = 0.435\textwidth]{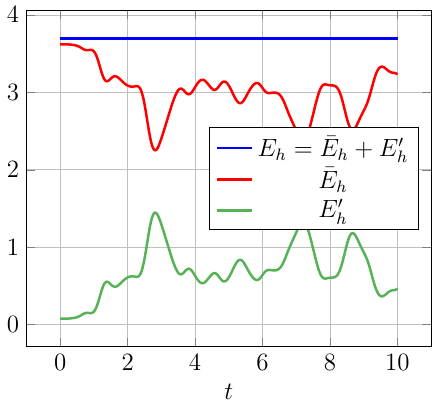}
    \caption{(Left) Fine-grid $\mathbf{u}$, reconstructed $\mathbf{R}\bar{\mathbf{u}}$, and SGS content $\mathbf{u}^\prime$ for $u = \sin(x)$. Here $N=1000$, $I = 20$, and $J=50$. The SGS content in the fourth coarse cell $\boldsymbol{\mu}_4$ is also indicated. (Right) Energy during a simulation of KdV equation with periodic BCs before and after filtering. }
    \label{fig:filtering}
\end{figure}

To reduce the degrees of freedom of the system we apply a spatial averaging filter to $\mathbf{u}$. This filter simply computes $\bar{\mathbf{u}}$ as a weighted average of $\mathbf{u}$, within each coarse cell. We represent this by the following matrix vector product:
\begin{equation}\label{eq:filtering}
    \bar{\mathbf{u}} = \mathbf{W}\mathbf{u},
\end{equation}
where $\mathbf{W}\in \mathbb{R}^{I\times N}$ is the filter.
The filter is defined as
\begin{equation}
    \mathbf{W} = \boldsymbol{\Omega}^{-1} \mathbf{O}  \boldsymbol{\omega}
\end{equation}
with overlap matrix $\mathbf{O}\in\mathbb{R}^{I\times N}$:
\begin{equation}
    \mathbf{O}:=
    \begin{bmatrix}
        1 & \hdots & 1 & & & & & &\\
        & & & \ddots &\ddots & \ddots & & & \\
        & & & & & & 1 & \hdots& 1
    \end{bmatrix}.
\end{equation}
This matrix contains ones at index $(i,j)$ if $\omega_{ij}$ lies in $\Omega_i$. 
Mathematically, $\bar{\mathbf{u}}$ can be regarded as a representation of $\mathbf{u}$ in a reduced basis. To project back onto the original basis we use a reconstruction operator $\mathbf{R}$:
\begin{equation}\label{eq:R}
    \mathbf{R} := \mathbf{O}^T
\end{equation}
which is a right inverse of $\mathbf{W}$, i.e.
\begin{equation}
    \mathbf{W}\mathbf{R} =\mathbf{I}.
\end{equation}
The matrix $\mathbf{R}$ approximates the reconstruction by a piece-wise constant function \cite{non_uniform_TRIAS2014246}. This is intuitively pictured in Figure \ref{fig:filtering}. 
Filtering a reconstructed solution $\mathbf{R} \bar{\mathbf{u}}$ leaves $\bar{\mathbf{u}}$ unchanged, i.e.
\begin{equation}\label{eq:proj_operator}
    \bar{\mathbf{u}} = \underbrace{(\mathbf{W}\mathbf{R})^p}_{=\mathbf{I}}\mathbf{W} \mathbf{u}
\end{equation}
for $p \in \mathbb{N}_0$. We will refer to this property as the `projection' property, as it is similar to how repeated application of a projection operator leaves a vector unchanged. 

By subtracting $\mathbf{R}\bar{\mathbf{u}}$ from $\mathbf{u}$ we obtain the subgrid-scale (SGS) content $\mathbf{u}^\prime \in \mathbb{R}^N$: \begin{equation}\label{eq:u_prime}
    \mathbf{u}^\prime := \mathbf{u} - \mathbf{R}\bar{\mathbf{u}}.
\end{equation}
In theory, one could define a more accurate $\mathbf{R}$, e.g.\ through polynomial reconstruction or data-driven approaches \cite{Syver}, and obtain a smaller $\mathbf{u}^\prime$.
However, this particular choice of $\mathbf{R}$ is made such that the energy is invariant under reconstruction, as will be shown in equation $\eqref{eq:filtering_integ_prop}$. This is an important property for our methodology to work.
Furthermore, we will refer to the SGS content in a single coarse cell $\Omega_i$ as $\boldsymbol{\mu}_i\in\mathbb{R}^{J}$, see Figure \ref{fig:filtering}.
Applying the filter to $\mathbf{u}^\prime$ yields zero:
\begin{equation}\label{eq:subgrid_zero}
    \mathbf{W}\mathbf{u}^\prime = \mathbf{W}\mathbf{u} - \underbrace{\mathbf{W}\mathbf{R}}_{=\mathbf{I}}\bar{\mathbf{u}}= \bar{\mathbf{u}} - \bar{\mathbf{u}} = 
    \mathbf{0}_\Omega,
\end{equation}
where $\mathbf{0}_\Omega$ is a vector with all entries equal to zero, defined on the coarse grid. This can be seen as the discrete equivalent of a property of a Reynolds operator \cite{sagaut2006large}. 
To illustrate, we display each of the introduced quantities for a 1D sinusoidal wave in Figure \ref{fig:filtering}.

\subsection{Discrete closure problem}

Next, we consider the time evolution of $\bar{\mathbf{u}}$. Since we employ a spatial filter which does not depend on time, filtering and time-differentiation commute: 
    $\mathbf{W}\frac{\text{d}\mathbf{u}}{\text{d} t} = \frac{\text{d} (\mathbf{W}\mathbf{u})}{\text{d} t}$. 
The closure problem arises because this is not true for the spatial discretization, i.e.
\begin{equation}
    \mathbf{W} f_h(\mathbf{u}) \neq f_H(\mathbf{W}\mathbf{u})
\end{equation}
where $f_H$ represents the same discretization scheme as $f_{h}$, but on the coarse grid. The closure problem is that the equations for $\bar{\mathbf{u}}$ are `unclosed'. This means that knowing the fine-grid solution $\mathbf{u}$ is required to evolve $\bar{\mathbf{u}}$ in time. In this way we do not achieve any computational speedup.

To resolve this we write the filtered system in closure model form:
\begin{align}
\label{eq:filtered_eqs}
    \frac{\text{d} \bar{\mathbf{u}}}{\text{d} t} = f_H(\bar{\mathbf{u}}) + \underbrace{(\mathbf{W}f_h(\mathbf{u})-f_H(\bar{\mathbf{u}}))}_{=: \mathbf{c}(\mathbf{u})},
\end{align}
where $\mathbf{c}(\mathbf{u})\in \mathbb{R}^I$ is the closure term. Note that this equation is still exact. $\mathbf{c}(\mathbf{u})$ is essentially the discrete equivalent of the commutator error in LES \cite{sagaut2006large}. One advantage of having first discretized the problem is that $\mathbf{c}(\mathbf{u})$ also includes the discretization error, with respect to a fine-grid simulation. The aim in closure modeling is generally to approximate $\mathbf{c}(\mathbf{u})$ by a closure model $\tilde{\mathbf{c}}(\bar{\mathbf{u}})$. In section \ref{sec:framework} we choose to represent $\tilde{\mathbf{c}}$ by a neural network.

\subsection{Inner products and energy decomposition} \label{sec:gen_E_content}

To describe the total energy that is present in the system, we define the following inner products and norms:
\begin{align}
     (\mathbf{a},\mathbf{b})_{\xi} &:= \mathbf{a}^T \boldsymbol{\xi}\mathbf{b} \\
    ||\mathbf{a}||^2_{\xi} &:= (\mathbf{a},\mathbf{a})_{\xi}
\end{align}
for $\boldsymbol{\xi}\in \{\boldsymbol{\omega},\boldsymbol{\Omega}\}$, and vectors $\mathbf{a}$ and $\mathbf{b}$.  With this notation we can represent the inner product on both the fine and coarse grid. For $\boldsymbol{\xi}=\mathbf{I}$ we obtain the conventional inner product and 2-norm, denoted as $(\mathbf{a},\mathbf{b}) = \mathbf{a}^T\mathbf{b}$ and $||\mathbf{a}||_2^2$. 
Besides the projection property \eqref{eq:proj_operator} an additional characteristic of the filter/reconstruction pair is that the inner product is conserved under reconstruction (see \ref{sec:filter_properties}):
\begin{equation}\label{eq:filtering_integ_prop}
    (\mathbf{R}\bar{\mathbf{a}},\mathbf{R}\bar{\mathbf{b}})_\omega = (\bar{\mathbf{a}},\bar{\mathbf{b}})_\Omega.
\end{equation}
Using this, the total energy $E_h$ in the system can be decomposed as
\begin{equation}\label{eq:E_decomp}
\begin{split}
    E_h :=& \frac{1}{2}||\mathbf{u}||^2_\omega =\frac{1}{2}||\mathbf{R}\bar{\mathbf{u}} + \mathbf{u}^\prime||^2_\omega 
    \\ =& \frac{1}{2}||\mathbf{R}\bar{\mathbf{u}}||^2_\omega + (\mathbf{R}\bar{\mathbf{u}}, \mathbf{u}^\prime)_\omega + \frac{1}{2}||\mathbf{u}^\prime||^2_\omega =\underbrace{ \frac{1}{2}||\bar{\mathbf{u}}||^2_\Omega}_{=:\bar{E}_h} + \underbrace{\frac{1}{2}|| \mathbf{u}^\prime||^2_\omega}_{=:E^\prime_h}
\end{split}
\end{equation}
where we replaced $\mathbf{u}$ by the decomposition in \eqref{eq:u_prime}. Furthermore, we used the fact that $\mathbf{R}\bar{\mathbf{u}}$ is orthogonal to $\mathbf{u}^\prime$ to simplify the expression, see \ref{sec:filter_properties}. The final expression shows that our choice of $\mathbf{W}$ and $\mathbf{R}$ is such that the total energy of the system can be split into two parts. These constitute of the resolved energy $\bar{E}_h$, which exclusively depends on $\bar{\mathbf{u}}$, and the SGS energy $E^\prime_h$, which exclusively depends on $\mathbf{u}^\prime$.
The energy conservation law can also be decomposed into a resolved and SGS part:
\begin{equation}\label{eq:local_energy_decomposition_evolution}
    \frac{\text{d}E_h}{\text{d}t} = \frac{\text{d}\bar{E}_h}{\text{d}t} + \frac{\text{d} E^\prime_h}{\text{d}t} = \left(\bar{\mathbf{u}},\frac{\text{d}\bar{\mathbf{u}}}{\text{d}t}\right)_\Omega + \left(\mathbf{u}^\prime,\frac{\text{d}\mathbf{u}^\prime}{\text{d}t}\right)_\omega =0,
\end{equation}
where we used the product rule to arrive at this relation. For Burgers' equation with $\nu>0$, the last equality sign changes to $\leq$. This means that even for dissipative systems the resolved energy can still increase (so-called `backscatter'), as long as the total energy is decreasing. For the KdV equation \eqref{eq:KdV}, which is strictly energy conserving, this decomposition can be seen in Figure \ref{fig:filtering}. Here one can clearly see the continuous exchange of energy between $\bar{E}_h$ and $E_h^\prime$, while the sum of the two remains constant.

\subsection{Momentum conservation}
Next to the energy, we investigate the effect of filtering on the momentum. The total discrete momentum is given by
\begin{equation}
    P_h = (\mathbf{1}_\omega,\mathbf{u})_\omega,
\end{equation}
where $\mathbf{1}_\omega$ is a vector with all entries equal to one, defined on the fine grid. From this definition we can show, see \ref{sec:filter_properties}, that the discrete momentum is invariant upon filtering, i.e.\
\begin{equation}\label{eq:mom_filter}
    (\mathbf{1}_\omega,\mathbf{u})_\omega  = (\mathbf{1}_\Omega,\bar{\mathbf{u}})_\Omega.
\end{equation}
This means the closure term does not add momentum into the system, i.e.
\begin{equation}\label{eq:mom_cons_closure}
(\mathbf{1}_\Omega,\mathbf{c}(\mathbf{u}))_{\Omega} = 0.
\end{equation}

\section{Structure-preserving closure modeling framework}\label{sec:framework}
In this section, the derived discrete energy and momentum balances will be used to construct a novel structure-preserving closure model. 
We will also discuss how to fit the parameters of the model. 

\subsection{The framework}
Many existing closure approaches aim at approximating $\mathbf{c}(\mathbf{u})$ by a closure model $\tilde{\mathbf{c}}(\bar{\mathbf{u}};\boldsymbol{\Theta})$. Here $\boldsymbol{\Theta}$ are parameters to be determined such that the approximation is accurate. In this work, we propose a novel formulation, in which we extend the system of equations for $\bar{\mathbf{u}}$ with a set of $I$ auxiliary SGS variables $\mathbf{s}\in \mathbb{R}^I$. These SGS variables locally approximate the SGS energy, but projected onto the coarse grid. This will be detailed later. The extended system of equations has the form
\begin{equation}\label{eq:full_eq}
    \frac{\text{d}}{\text{d}t}\begin{bmatrix}
    \bar{\mathbf{u}} \\ \mathbf{s}
    \end{bmatrix}\approx \mathcal{G}_{\Theta}(\bar{\mathbf{u}},\mathbf{s}) := \begin{bmatrix}
    f_H(\bar{\mathbf{u}}) \\
    \mathbf{0}_\Omega
    \end{bmatrix} + \boldsymbol{\Omega}_2^{-1}(\mathcal{K}-\mathcal{K}^T)\begin{bmatrix}
    \bar{\mathbf{u}} \\ \mathbf{s}
    \end{bmatrix} - \boldsymbol{\Omega}_2^{-1}\mathcal{Q}^T
    \mathcal{Q} \begin{bmatrix}
    \bar{\mathbf{u}} \\ \mathbf{s}
    \end{bmatrix},
\end{equation}
where $\mathcal{K}=\mathcal{K}(\bar{\mathbf{u}},\mathbf{s},\mathbf{\Theta})\in \mathbb{R}^{2I\times 2I}$ and $\mathcal{Q}=\mathcal{Q}(\bar{\mathbf{u}},\mathbf{s},\mathbf{\Theta})\in \mathbb{R}^{2I\times 2I}$ depend on the solution in a parameterized fashion.
Next to the introduction of $\mathbf{s}$, the second main novelty in this work is to formulate the closure model in terms of a skew-symmetric and a dissipative term. The skew-symmetric term is introduced to allow for exchange of energy between $\bar{\mathbf{u}}$ and $\mathbf{s}$. The dissipative term is introduced to provide additional dissipation, as this is required (see \ref{sec:diss-diff}). The operators $\mathcal{K}$ and $\mathcal{Q}$ will be modeled in terms of neural networks (NNs) with trainable parameters (contained in $\boldsymbol{\Theta}$). So even though the notation in \eqref{eq:full_eq} suggests linearity of the closure model, the dependence of $\mathcal{K}$ and $\mathcal{Q}$ on $\bar{\mathbf{u}}$ and $\mathbf{s}$ makes the model nonlinear. The construction of the introduced operators will be detailed in sections \ref{sec:defining_operators}. As our energy definition includes $\boldsymbol{\Omega}$ we include the inverse of the concatenated mass matrix $\boldsymbol{\Omega}^{-1}_2$ in our system of equations \eqref{eq:full_eq}. This ensures energy conservation/dissipation regardless of the grid topology. This mass matrix is defined as 
\begin{equation}\label{eq:Omega_2}
    \boldsymbol{\Omega}_2 = \begin{bmatrix}
\boldsymbol{\Omega} & \\ & \boldsymbol{\Omega}
    \end{bmatrix}.
\end{equation}

Given the extended system of equations, the total energy \eqref{eq:E_decomp} is approximated as 
\begin{equation}\label{eq:E_s}
    E_h \approx E_s := \frac{1}{2}||\mathbf{a}||^2_{\Omega_2} = \frac{1}{2}||\bar{\mathbf{u}}||^2_\Omega+ \frac{1}{2}||\mathbf{s}||^2_\Omega,
\end{equation}
where the second term approximates the SGS energy.
Furthermore, we concatenate $\bar{\mathbf{u}}$ and $\mathbf{s}$ into a single state vector $\mathbf{a} \in \mathbb{R}^{2I}$:
\begin{equation}\label{eq:bar_U}
    \mathbf{a} := \begin{bmatrix}
    \bar{\mathbf{u}} \\ \mathbf{s}
    \end{bmatrix}.
\end{equation}
The evolution equation for the approximated total energy is given by
\begin{equation}\label{eq:en_cons_closure}
    \frac{\text{d}E_s}{\text{d}t} 
    = \left(\mathbf{a},\frac{\text{d}\mathbf{a}}{\text{d}t}    \right)_{\Omega_2} = (\bar{\mathbf{u}},f_H(\bar{\mathbf{u}}))_\Omega - ||\mathcal{Q}\mathbf{a}||_2^2  ,
\end{equation}
as the skew-symmetric term involving $\mathcal{K} - \mathcal{K}^T$ cancels. This is a property of skew-symmetric matrices \cite{benjamin_thesis}. Consequently, this formulation guarantees stability provided that $f_{H}$ is structure-preserving.

Our key insight is that \textit{by explicitly including an approximation of the SGS energy we are able to satisfy the energy conservation balance, equation \eqref{eq:local_energy_decomposition_evolution}}. The energy balance serves not only as an important constraint for the closure model (represented by a NN), but also guarantees stability of our closure model. This is because the energy is a norm of the solution which is bounded in time. This framework thus allows for the modeling of backscatter without sacrificing stability.

\subsection{SGS variables}\label{sec:SGS_variables}

Next, let us consider appropriate expressions for $\mathbf{s}$. The exact SGS energy on the coarse grid is given by:
\begin{equation}\label{eq:SGS_energy}
    \mathbf{W}(\mathbf{u}^{\prime})^2,
\end{equation}
where $(.)^2$ is to be interpreted element-wise. This would yield $\mathbf{s}= \pm\sqrt{\mathbf{W}(\mathbf{u}^\prime)^2}$, where $\sqrt{(.)}$ is also to be interpreted element-wise. The square root is taken to comply with the energy definition in \eqref{eq:E_s}. However, this definition for $\mathbf{s}$ resulted in poor performance during testing. We argue this was caused by the strict positivity (or negativity). Inevitable small errors in the model predictions caused some of the elements of $\mathbf{s}$ to switch sign. As the model was trained on positive $\mathbf{s}$, the simulation quickly diverged from the true trajectory. Attempts at resolving this issue were not successful.

Instead we propose the use of a local linear compression. This formulation naturally allows for both positive and negative values of $\text{s}_i$.
This compression is written as (assuming a uniform grid):
\begin{equation}\label{eq:def_s}
    \text{s}_i = \mathbf{t}^T\boldsymbol{\mu}_i, \qquad i=1,\ldots,I,
\end{equation}
where we recall that $\boldsymbol{\mu}_i\in\mathbb{R}^{J}$ represents the SGS content in a single coarse cell $\Omega_i$.
Furthermore, $\mathbf{t}\in \mathbb{R}^J$ are the compression parameters. 
We aim to choose $\mathbf{t}$ such that we obtain $\mathbf{s}^2 \approx \mathbf{W}(\mathbf{u}^\prime)^2$. The optimal values of $\mathbf{t}$ are obtained using a singular value decomposition of the SGS content. This is outlined in \ref{sec:SGS_compression}.
From \eqref{eq:def_s} we construct an operator $\mathbf{T}_s\in \mathbb{R}^{I \times N}$ which transform the SGS content into $\mathbf{s}$:
\begin{equation}
    \mathbf{s} = \mathbf{T}_s \mathbf{u}^\prime. 
\end{equation}
Combining the compression with the filter, see \eqref{eq:filtering}, we define the operator $\mathbf{T}\in \mathbb{R}^{2I \times N}$ which transforms $\mathbf{u}$ into the state vector $\mathbf{a}$:
\begin{equation}
    \mathbf{a}= \underbrace{\begin{bmatrix}
        \mathbf{W} \\
        \mathbf{T}_s(\mathbf{I}-\mathbf{R}\mathbf{W})
    \end{bmatrix}}_{=:\mathbf{T}}\mathbf{u}.
\end{equation}
Due to the linearity of the transformation we simply obtain
\begin{equation}
    \frac{\text{d}\mathbf{a}}{\text{d}t}= \mathbf{T}\frac{\text{d}\mathbf{u}}{\text{d}t},
\end{equation}
where $\mathbf{T}$ is the Jacobian of the transformation. In the linear case this Jacobian does not depend on $\mathbf{u}$ which significantly simplifies computing reference data for $\frac{\text{d}\mathbf{a}}{\text{d}t}$.
In addition, it follows that if the true RHS includes a forcing term $\mathbf{F} \in \mathbb{R}^N$ we simply account for this by adding $\mathbf{T}\mathbf{F}$
to the RHS of \eqref{eq:full_eq}.

\begin{figure}
    \centering
    \begin{subfigure}[b]{0.441\textwidth}
        \includegraphics[width = \textwidth]{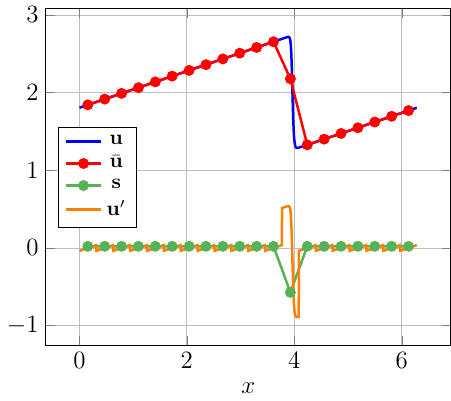}
    \end{subfigure}
    \begin{subfigure}[b]{0.48\textwidth}
        \includegraphics[width = \textwidth]{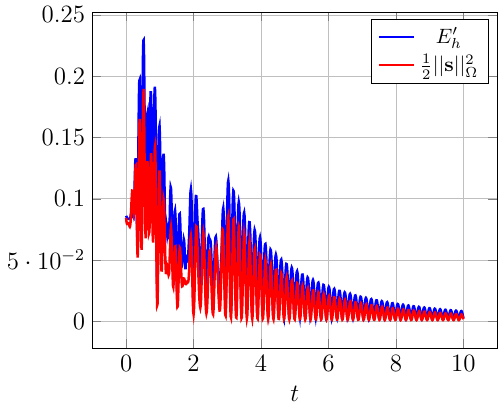}
    \end{subfigure}
    \caption{(Left) Learned SGS compression applied to Burgers' equation for $N=1000$, with $I=20$ and $J=50$. By filtering and applying the SGS compression the degrees of freedom of this system are effectively reduced from $N=1000$ to $2I = 40$. (Right) True SGS energy and compressed SGS energy during this simulation of Burgers' equation.}
    \label{fig:compression}
\end{figure}

To illustrate how the compression works in practice we consider a snapshot from a simulation of Burgers' equation ($\nu = 0.01$) with periodic BCs, see Figure \ref{fig:compression}.
We observe that $\mathbf{s}$ serves as an energy storage for the SGS content, which is mainly present near shocks. If we look at the SGS energy trajectory we find that its behavior is captured both qualitatively and quantitatively by the SGS compression. Although we still miss some of the energy, the oscillations due to the traveling shock are nicely captured. From this we argue that a linear compression suffices. For more complex systems autoencoders might offer an alternative \cite{bank2023autoencoders}.

\subsection{Construction of the operators}\label{sec:defining_operators}

Similarly to the structure-preserving discretization for Burgers' equation, presented in section \ref{sec:burgers}, we want our closure model to locally advect momentum and energy through the domain. In this way we do not violate the conservation laws. It is therefore that we inspire our machine learning closure model on this structure-preserving discretization. In this section we outline how to construct $\mathcal{K}$ and $\mathcal{Q}$ from the output of a convolutional neural network (CNN) \cite{conv_NN}.

\subsubsection{Diffusive operator}

Let us start by considering the diffusion operator in \eqref{eq:Q^TQ}, namely $-\mathbf{Q}^T\text{diag}(\boldsymbol{\nu})\mathbf{Q}$. This operator locally diffuses energy and momentum through space at a rate which is determined by a positive viscosity field. In the periodic case the momentum is conserved, as $\mathbf{1}^T \mathbf{Q}^T = \mathbf{0}$. and the energy is dissipated, as $-\mathbf{u}^T\mathbf{Q}\text{diag}(\boldsymbol{\nu})\mathbf{Q}\mathbf{u}= -||\sqrt{\text{diag}(\boldsymbol{\nu})}\mathbf{Q}\mathbf{u}||_2^2 \leq 0$.

Drawing inspiration from this operator we introduce the following form for the diffusive term in our closure model:
\begin{equation}\label{eq:Q}
     \mathcal{Q}(\mathbf{a};\boldsymbol{\Theta}) = \mathbf{q}(\mathbf{a};\boldsymbol{\Theta})\mathcal{B}_1(\boldsymbol{\Theta})
\end{equation}
such that $\mathcal{Q}^T\mathcal{Q} = \mathcal{B}_1^T\mathbf{q}^2\mathcal{B}_1$ resembles the discussed diffusion operator. Here $\mathbf{q}(\mathbf{a};\boldsymbol{\Theta}) = \text{diag}(\mathbf{q}_1,\mathbf{q}_2)\in \mathbb{R}^{2I\times 2I}$ is constructed from two output channels $\mathbf{q}_1$  and $\mathbf{q}_2$ of a CNN which takes $\bar{\mathbf{u}}, \mathbf{s}$, and $f_H(\bar{\mathbf{u}})$ as inputs. $f_H(\bar{\mathbf{u}})$ was added as input channel because it significantly improved the performance of the neural network. This is supported by \cite{Melchers2022}. Looking at the introduced form, $\mathbf{q}^2$ can be thought of as a set of learned non-uniform and nonlinear viscosity fields. The square ensures positivity of these fields. Furthermore, the introduced matrix $\mathcal{B}_1(\boldsymbol{\Theta})\in \mathbb{R}^{2I\times 2I}$ is a linear operator which encodes a set of parameterized convolutions. It will be used to satisfy momentum conservation, similarly to $\mathbf{Q}$. This will be discussed later in section \ref{sec:momentum_cons}.

\subsubsection{Advective operator}

For our skew-symmetric operator $\mathcal{K} - \mathcal{K}^T$ we take a similar approach and introduce the following form:
\begin{equation}
    \mathcal{K}(\mathbf{a};\boldsymbol{\Theta})= \mathcal{B}^T_2(\boldsymbol{\Theta})\mathbf{k}(\mathbf{a};\boldsymbol{\Theta})\mathcal{B}_3(\boldsymbol{\Theta})
\end{equation}
such that $\mathcal{K} - \mathcal{K}^T = \mathcal{B}^T_2\mathbf{k}\mathcal{B}_3 - \mathcal{B}^T_3\mathbf{k}\mathcal{B}_2$. Similarly to the diffusive operator we use $\mathcal{B}_2(\boldsymbol{\Theta}),\mathcal{B}_3(\boldsymbol{\Theta}) \in \mathbb{R}^{2I\times2I}$ to satisfy momentum conservation. As is the case for $\mathbf{q}$, the fields $\mathbf{k}(\mathbf{a};\boldsymbol{\Theta}) = \text{diag}(\mathbf{k}_1,\mathbf{k}_2)$ are constructed from an additional set of two outputs channels of the CNN, i.e.
\begin{equation}
    \begin{bmatrix}
        \mathbf{q}_1 & \mathbf{q}_2 & \mathbf{k}_1 & \mathbf{k}_2
    \end{bmatrix} = \text{CNN}(\bar{\mathbf{u}},\mathbf{s},f(\bar{\mathbf{u}});\boldsymbol{\Theta}).
\end{equation}
This means the CNN has in total four output channels to construct both $\mathbf{q}$ and $\mathbf{k}$.
Furthermore, $\mathbf{k}$ can be thought of as a set of learned velocity fields which advect momentum and energy through the domain, in addition to exchanging energy between $\bar{\mathbf{u}}$ and $\mathbf{s}$.

\subsubsection{Momentum conservation}\label{sec:momentum_cons}

The entire framework \eqref{eq:full_eq} is summarized as
\begin{equation}\label{eq:full_eq_plus}
     \mathcal{G}_{\Theta}(\mathbf{a}) = \begin{bmatrix}
    f_H(\bar{\mathbf{u}}) \\
    \mathbf{0}_\Omega
    \end{bmatrix} + \boldsymbol{\Omega}_2^{-1}(\mathcal{B}^T_2\mathbf{k}\mathcal{B}_3 - \mathcal{B}^T_3\mathbf{k}\mathcal{B}_2)\mathbf{a} - \boldsymbol{\Omega}_2^{-1}\mathcal{B}_1^T\mathbf{q}^2\mathcal{B}_1 \mathbf{a}.
\end{equation}
From this we find that momentum conservation, see \eqref{eq:mom_cons_closure}, places an additional constraint on the operators:
\begin{equation}\label{eq:skew_closure_mom_cons}
\left(\begin{bmatrix}
    \mathbf{1}_\Omega \\ \mathbf{0}_\Omega
\end{bmatrix},\mathcal{G}_\Theta(\mathbf{a})\right)_{\Omega_2} = \begin{bmatrix}
    \mathbf{1}_\Omega \\ \mathbf{0}_\Omega
\end{bmatrix}^T(\mathcal{B}^T_2\mathbf{k}\mathcal{B}_3 - \mathcal{B}^T_3\mathbf{k}\mathcal{B}_2)\mathbf{a} - \begin{bmatrix}
    \mathbf{1}_\Omega \\ \mathbf{0}_\Omega
\end{bmatrix}^T\mathcal{B}_1^T\mathbf{q}^2\mathcal{B}_1 \mathbf{a} = 0,
\end{equation}
assuming that $f_{H}$ is momentum conserving. This constraint is to be satisfied for periodic BCs, and we choose to compose $\mathcal{B}_i$ of convolutions (or stencils) expressed as linear operators:
\begin{equation}\label{eq:B_matrix}
    \mathcal{B}_i= \begin{bmatrix}
        \mathbf{B}_i^{11} & \mathbf{B}_i^{12}\\ \mathbf{B}_i^{21} & \mathbf{B}_i^{22}
    \end{bmatrix}.
\end{equation} 
The operator $\mathcal{B}_i$ can therefore be thought of as the connection between two layers in a CNN, with each layer containing two channels \cite{conv_NN}. In 1D each of the submatrices is characterized by $2B+1$ parameters, where $B>0$ is the width of the convolution. In this way, each convolution takes into account $B$ neighboring grid cells from each side. Momentum conservation is ensured by constraining these submatrices in a clever way. From the constraint \eqref{eq:skew_closure_mom_cons} and the definition of $\mathcal{B}_i$, in \eqref{eq:B_matrix}, we find that momentum conservation is satisfied if the sum of the convolution weights for $\mathbf{B}^{j1}_i$ is zero $\forall i,j$, as is the case for $\mathbf{Q}$. To see how this works we consider a general convolution operator $\mathbf{B}$, characterized by parameters $b_{-B},\ldots,b_{B} \in \mathbb{R}$. Applying this operator to a discrete field $\mathbf{f}$, while constraining the sum of the weights to zero, is achieved as follows:
\begin{align}
    (\mathbf{B}\mathbf{f})_i &= \sum^B_{j=-B} \bar{b}_j\text{f}_{i+j}, \\
    \bar{b}_j &= b_j - \sum^B_{k=-B}\frac{b_k}{2B+1},
\end{align}
such that $\sum^{B}_{j=-B}\bar{b}_j = 0$ indeed holds.
Applying this procedure to $\mathbf{B}^{j1}_i$, $\forall i,j$, ensures \eqref{eq:skew_closure_mom_cons} is satisfied up to machine precision, independent of the parameter values. The remaining convolutions are left unconstrained, i.e.\ we simply take $\bar{b}_j = b_j$.

\subsubsection{Properties \& further discussion}

The key insight is that we are free to choose any set of parameters $\boldsymbol{\Theta}$ without violating the prescribed structure of the system. Furthermore, as $\mathcal{K}$ and $\mathcal{Q}$ are based solely on convolutions in the CNN and the $\mathcal{B}$ matrices, they effectively correspond to nonlinear local stencils. The entire framework is therefore translation equivariant. For periodic boundary conditions (BCs) we apply circular padding to both the CNN inputs and the state vector \cite{conv_NN}. For non-periodic BCs we refer to \ref{sec:BCs}. The parameters $\boldsymbol{\Theta}$ include the weights of the CNN, as well as the parameters characterizing the $\mathcal{B}$ matrices. As both the CNN and the convolution operations in $\mathcal{B}$ are sparse, our model remains computationally efficient. 

\begin{figure}
    \centering
    \begin{subfigure}[b]{0.9 \textwidth}
        \includegraphics[width = \textwidth]{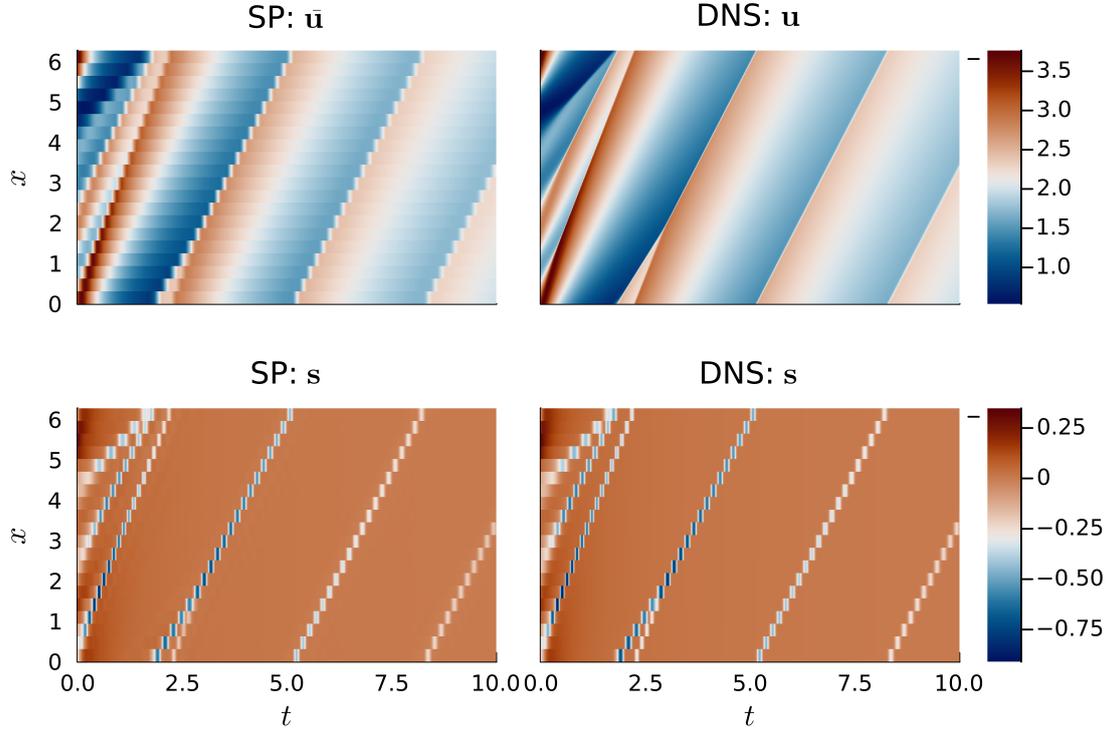}
    \end{subfigure}
    \caption{A simulation of Burgers' equation with periodic BCs using our trained structure-preserving closure model for $I = 20$ and $J = 50$ (left), along with the DNS solution for $N=1000$ (right).}
    \label{fig:heatmap_per}
\end{figure}

\begin{figure}
    \centering
    \begin{subfigure}[b]{0.50\textwidth}
        \includegraphics[width = \textwidth]{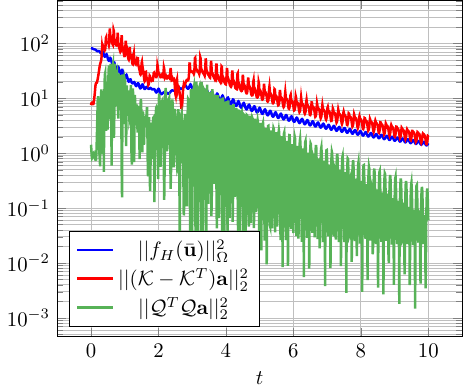}
    \end{subfigure}
    \caption{Magnitude of each of the different terms present in \eqref{eq:full_eq} corresponding to the simulation in Figure \ref{fig:heatmap_per} with $I = 20$, $J = 50$, and $N=1000$.}
    \label{fig:decomp}
\end{figure}

In Figure \ref{fig:heatmap_per} the framework is applied to Burgers' equation, where we compare it to the direct numerical simulation (DNS). 
It is once again interesting to see that $\mathbf{s}$ is largest at the shocks, indicating the presence of significant SGS content. Comparing the magnitude of the different terms in \eqref{eq:full_eq}, see Figure \ref{fig:decomp}, we observe that the skew-symmetric term, that is responsible for redistributing the energy, is most important. In fact it is more important than the coarse-grid discretization. In other words, our closure model has learned dynamics that are highly significant to correctly predict the evolution of the filtered system. This means that even though the closure term is large, we can still accurately model it.

\subsection{Finding the optimal parameter values}

The optimal parameter values of the network can be obtained numerically by minimizing 
\begin{equation}\label{eq:derivative_loss}
\mathcal{L}(\mathbf{X}_\text{u};\boldsymbol{\Theta}):=\frac{1}{p}\sum_{\mathbf{u} \in \mathbf{X}_\text{u}} ||\mathcal{G}_{\Theta}(\mathbf{T}\mathbf{u})-\mathbf{T}f_h(\mathbf{u})||_{2}^2
\end{equation}
with respect to $\boldsymbol{\Theta}$ for the training set $\mathbf{X}_\text{u}$ containing $p$ DNS snapshots. 
We will refer to this approach as `derivative fitting', as we minimize the residual between the predicted and the true RHS.

An alternative is to optimize $\boldsymbol{\Theta}$ such that the solution itself is accurately reproduced. To achieve this we minimize
\begin{equation}\label{eq:trajectory_loss}
\mathcal{L}_n(\mathbf{X}_\text{u};\boldsymbol{\Theta}):=\frac{1}{pn}\sum_{\mathbf{u}\in\mathbf{X}_\text{u}}\sum_{i=1}^n || \bar{\mathcal{S}}^i_{\boldsymbol{\Theta}}(\mathbf{T}\mathbf{u})  - \mathbf{T}\mathcal{S}^{i(\overline{\Delta t}/\Delta t) }(\mathbf{u}) ||_{2}^2.
\end{equation}
Here $\bar{\mathcal{S}}_{\boldsymbol{\Theta}}^i(\mathbf{T}\mathbf{u})$ represents the output of the solver after successive application of an explicit time integration scheme for $i$ time steps, with step size $\overline{\Delta t}$, starting from initial condition $\mathbf{T}\mathbf{u}$, using the introduced closure model in \eqref{eq:full_eq}. The DNS counterpart is indicated by $\mathcal{S}^{i(\overline{\Delta t}/\Delta t)}(\mathbf{u})$, with step size $\Delta t$, starting from initial condition $\mathbf{u}$.  Note the appearance of the ratio $\overline{\Delta t}/\Delta t$. This is because the coarse grid allows for larger time steps \cite{de2013courant_CFL}.  We will refer to this method as `trajectory fitting'. This approach has been shown to yield more accurate and stable closure models \cite{Frezat,List,MacArt_embedded_learning,Syver,Melchers2022}. 
In this paper we propose a hybrid of the two approaches, as trajectory fitting is more computationally more expensive. This hybrid approach will be detailed later.

\section{Results}\label{sec:results}

\subsection{Test cases}

To test our closure modeling framework we consider the previously introduced Burgers' equation with $\nu = 0.01$  on the spatial domain $\Omega = [0,2\pi]$ for two test cases: (i) periodic BCs without forcing and (ii) inflow/outflow (I/O) BCs with time-independent forcing. The implementation of BCs is discussed in \ref{sec:BCs}. We also consider a third test case: (iii) the KdV equation with $\varepsilon = 6$ and $\mu = 1$ on the spatial domain $\Omega = [0,32]$ for periodic BCs.
Parameter values for Burgers' and KdV are taken from \cite{Bar-Sinai}. Reference DNSs are carried out on a uniform grid of $N=1000$ for Burgers' and $N=600$ for KdV up to time $T=10$. The data that is generated from these reference simulations is split into a training (70\%) and a validation set (30\%). The simulation conditions (initial conditions, BCs, and forcing) for training are generated randomly, as described in \ref{sec:training}.
The closure models will be tested on unseen simulation conditions, sampled from the same distributions as the training data. 
In addition to this, the construction of the training and validation set, the training procedure, and the hyperparameter tuning procedure are also described in \ref{sec:training}.

\subsubsection{Considered closure models}

For the analysis, we will compare our structure-preserving framework (SP) to a vanilla CNN. The vanilla CNN output is multiplied by a forward difference operator and the result is used as a closure model, i.e.
\begin{equation}\label{eq:CNN_closure}
    \tilde{\mathbf{c}}(\mathbf{u};\boldsymbol{\Theta}) = \bar{\mathbf{Q}}\text{CNN}(\bar{\mathbf{u}},f_h(\bar{\mathbf{u}});\boldsymbol{\Theta}),
\end{equation}
where $\bar{\mathbf{Q}}$ is a forward difference discretization of the first derivative on the coarse grid. This is done to satisfy momentum conservation \cite{Melchers2022,burgers_3}. In addition, we consider a constant Smagorinsky model (SM) for use case (i) and (ii). This model contains a single scalar parameter $C_s$. This parameter is optimized in the same way as the NN parameters. We chose this model as it was the best performing non-machine learning closure in \cite{burgers_3}, outperforming the dynamic Smagorinsky model.
After discretization, it takes the following form:
\begin{align}
    \tilde{\mathbf{c}}(\mathbf{u};C_s) &= -\bar{\mathbf{Q}}^T\text{diag}(\boldsymbol{\nu}_t(C_s))\bar{\mathbf{Q}}\bar{\mathbf{u}}, \\
    \boldsymbol{\nu}_t(C_s) &= (hC_s)^2 |\bar{\mathbf{Q}}\bar{\mathbf{u}}|,
\end{align}
such that $\boldsymbol{\nu}_t$ represents a parameterized and solution dependent viscosity. Here we have taken the filter width equal to the grid-spacing $h$. Moreover, we also consider SP0 which initializes the simulation with $\mathbf{s}(t=0)= \mathbf{0}_\Omega$ instead of the true $\mathbf{s}$. This emulates the situation in which the true initial condition is unknown.
Finally, we consider the no closure (NC) case, which corresponds to a coarse-grid solution of the PDEs. 
To march the solution forward in time we employ an explicit RK4 scheme \cite{RK4_Butcher:2007} with time step size $\overline{\Delta t}=0.01$ ($4\times$ larger than the DNS) for Burgers' and $\overline{\Delta t}=5 \times 10^{-3}$ ($50\times$ larger than the DNS) for KdV. 

To make a fair comparison we compare closure models with the same number of degrees of freedom (DOF). For SP we have $\text{DOF}=2I$, as we obtain an additional set of $I$ degrees of freedom corresponding to the addition of the SGS variables. For the remaining closure models we simply have $\text{DOF}=I$. 

\subsection{Training the closure models}

As use cases (i) and (ii) both correspond to Burgers' equation we train the corresponding closure models on a dataset containing both simulation conditions. In this way we end up with a single closure model that works for both (i) and (ii). The SP closure models contain in total 2780 parameters (consisting of two hidden layers with each 20 channels and a kernel size of 5 for the underlying CNN) for Burgers' equation and 5352 (consisting of two hidden layers with each 30 channels and a kernel size of 5) for KdV. The purely CNN-based closure models consist of 3261 parameters (two hidden layers with each 20 channels and a kernel size of 7). These settings are based on the hyperparameter tuning procedure in \ref{sec:training}. For KdV we omit the dissipative component in \eqref{eq:full_eq}, as it is a conservative system. In between hidden layers we employ the ReLU activation function. We employ a linear activation function at the final layer. For Burgers' we choose $B=1$ for the construction of the $\mathcal{B}$ matrices, matching the width of the coarse discretization. For KdV we do the same and therefore take $B=2$. 

As stated earlier, the model parameters are optimized by first derivative fitting and then trajectory fitting. During testing we observed that solely derivative fitting resulted in instabilities and poor performance for the vanilla CNN, especially in the KdV case. In contrast, our SP method is guaranteed to be stable, and simply derivative fitting already resulted in reasonable performance. To maximize the potential of each closure model we decided to include trajectory fitting in the training procedure. The full procedure is outlined in \ref{sec:training}. We implemented our closure models in the Julia programming language \cite{bezanson2017julia} using the Flux.jl package \cite{Flux.jl-2018,flux_innes:2018}. The code can be found at \url{https://github.com/tobyvg/ECNCM_1D}.

\subsection{Closure model performance}

We examine the performance of the trained closure models based on how well the filtered DNS solution is reproduced. For our comparison we will make extensive use of the normalized root-mean-squared error (NRMSE) metric. This metric is defined as
\begin{equation}
    \text{NRMSE } \bar{\mathbf{u}} (t)=\sqrt{\frac{1}{|\Omega|}||\bar{\mathbf{u}}(t)-\bar{\mathbf{u}}^{\text{DNS}}(t)||^2_\Omega}.
\end{equation}
It is used to compare the approximated solution $\bar{\mathbf{u}}$ at time $t$, living on the coarse grid, to the filtered DNS result $\bar{\mathbf{u}}^\text{DNS}$. We will refer to this metric as the solution error.
 In addition, we define the integrated-NRMSE 
 (I-NRMSE) as
 \begin{equation}\label{eq:int_sol}
     \text{I-NRMSE }\bar{\mathbf{u}} = \frac{1}{T} \sum_{i}\overline{\Delta t} \text{ NMRSE } \bar{\mathbf{u}} (i\overline{\Delta t}), \qquad 0\leq i\overline{\Delta t} \leq T,
 \end{equation}
such that the sum represents integrating the solution error in time. We will refer to this metric as the integrated solution error.

\subsubsection{Convergence}

As we refine the resolution of the coarse grid, and with this increase the number of $\text{DOF}$, we expect convergence of both the compression error 
$\mathcal{L}_s$, see \ref{sec:SGS_compression}, and the solution error. We consider $\text{DOF}\in \{20,30,40,50,60,70,80,90,100\}$. For each $
\text{DOF}$ value we train a different set of closure models.
If $N$ is not divisible by $I$ we first project the DNS result onto a grid with a resolution that is divisible by $I$ before applying the filter.
In total 36 closure models are trained: two (SP and CNN) for each combination of the 9 considered coarse-grid resolutions and equation (Burgers' and KdV). 

\begin{figure}
    \centering
    \begin{subfigure}[b]{0.414\textwidth}
        \includegraphics[width = \textwidth]{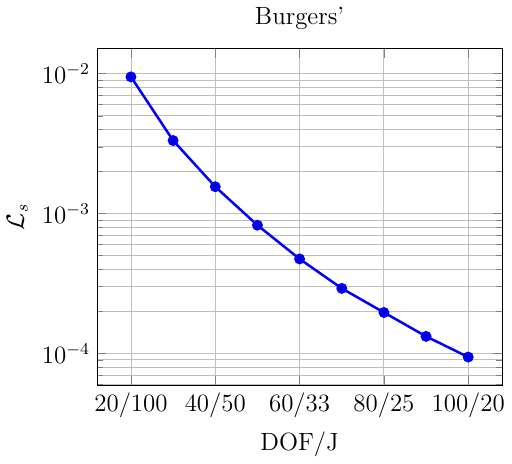}
    \end{subfigure}
    \begin{subfigure}[b]{0.386\textwidth}
        \includegraphics[width = \textwidth]{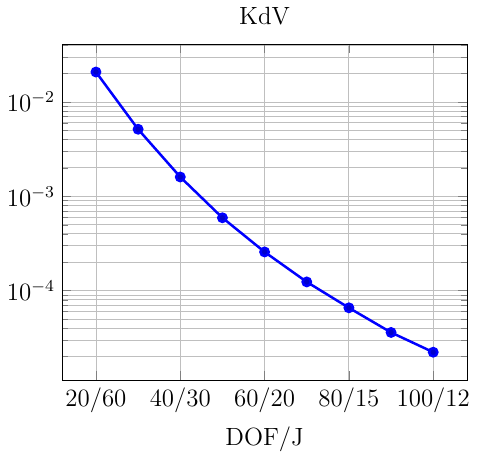}
    \end{subfigure}
    \caption{Convergence of the SGS compression error $\mathcal{L}_s$ when refining the coarse grid, evaluated on the validation set for Burgers' equation ($N=1000$) and KdV equation ($N=600$). Both the effective $\text{DOF}$ and the corresponding compression factor $J$ are depicted on the $x$-axis.}
    \label{fig:compression_convergence}
\end{figure}

The SGS compression error evaluated over the validation set is shown in Figure \ref{fig:compression_convergence}.
We observe monotonic convergence of the compression error as we refine the grid. We expect the compression error to further converge to zero, as we keep on refining. The faster convergence for the KdV equation is likely caused by the lower fine-grid resolution, as compared to Burgers' equation. 

\begin{figure}[ht]
    \centering
    \begin{subfigure}[b]{0.326\textwidth}
        \includegraphics[width = \textwidth]{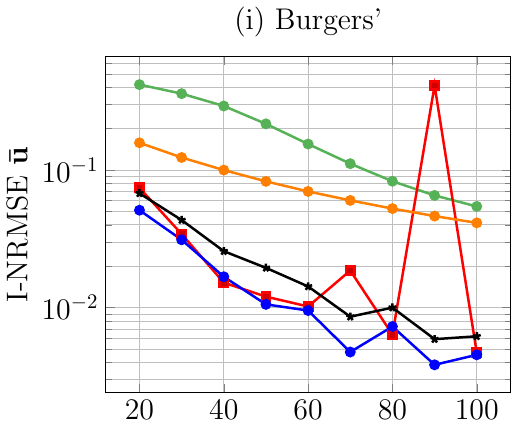}
    \end{subfigure}
    \begin{subfigure}[b]{0.3045\textwidth}
        \includegraphics[width = \textwidth]{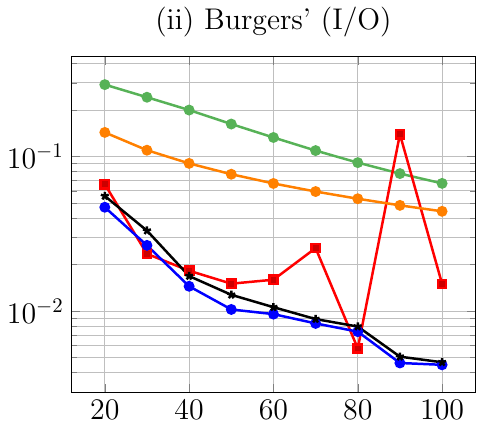}
    \end{subfigure}
    \begin{subfigure}[b]{0.3045\textwidth}
        \includegraphics[width = \textwidth]{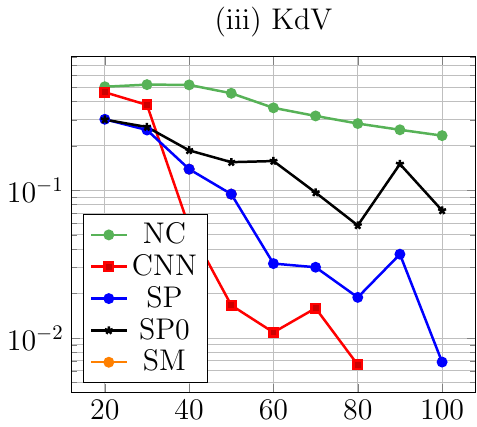}
    \end{subfigure}
    \begin{subfigure}[b]{0.33\textwidth}
        \includegraphics[width = \textwidth]{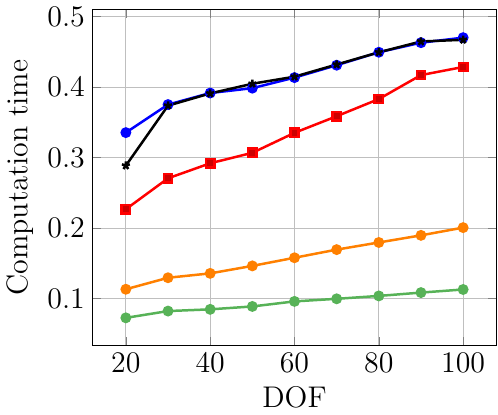}
    \end{subfigure}
    \begin{subfigure}[b]{0.304\textwidth}
        \includegraphics[width = \textwidth]{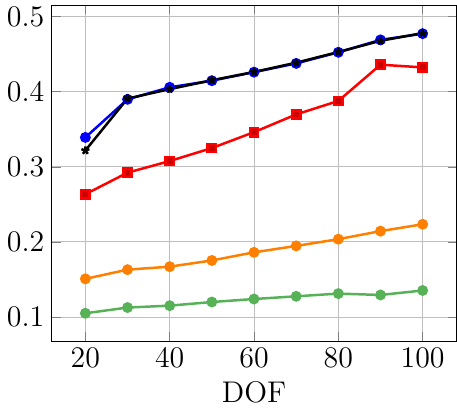}
    \end{subfigure}
    \begin{subfigure}[b]{0.289\textwidth}
        \includegraphics[width = \textwidth]{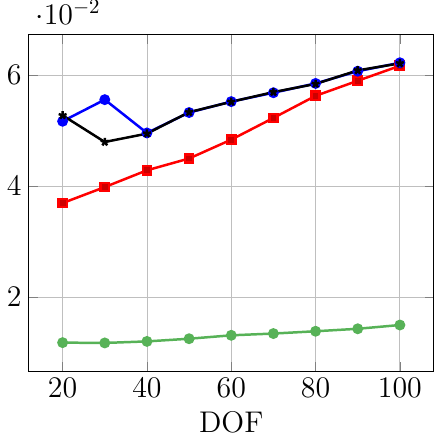}
    \end{subfigure}
    \caption{(Top) Integrated solution error evaluated at $T=10$ averaged over 20 simulations for the different use cases (i)-(iii) and an increasing number of $\text{DOF}$. Only stable simulations are considered for the depicted averages. Absence of a scatter point indicates none of the simulations were stable. (Bottom) Similar to the top plots, but depicting the average simulation time as a fraction of the DNS time. $\textbf{NC}=\text{no closure}$, $\textbf{CNN}=\text{convolutional neural network closure}$, $\textbf{SP}=\text{structure-preserving closure}$, $\textbf{SP0}=\text{structure-preserving closure with $\mathbf{s}(t=0)=\mathbf{0}_\Omega$}$, and $\textbf{SM}=\text{constant Smagorinsky model}$.}
    \label{fig:convergence}
\end{figure}

Next, we consider the integrated solution error, see Figure \ref{fig:convergence}. The presented plots represent an average taken over 20 simulations with unseen simulation conditions. We consider different numbers of $\text{DOF}$ which enables us to evaluate the convergence.
For test cases (i) and (ii) we observe almost monotonic convergence of the solution error for NC, SM, and SP/SP0.
Furthermore, SP improves upon NC with roughly one order of magnitude, surpassing SM. SP0 only performs slightly worse than SP, but still converges relatively smoothly. On the other hand, the solution error for the CNN behaves quite erratically: sometimes more accurate than SP, and sometimes less accurate than NC (case (i), $\text{DOF}=90$). 

For test case (iii) (KdV) we find that for most numbers of $\text{DOF}$ the CNN outperforms SP, while not resulting in stable closure models for $\text{DOF}\in \{90,100\}$. Furthermore, for the lower numbers of $\text{DOF}$ we observe slightly better performance for SP. From this we conclude that the compression error (see Figure \ref{fig:compression_convergence}) is likely not the limiting factor of the closure model performance. 
Looking at the difference between SP and SP0 it seems that initializing with the true $\mathbf{s}$ seems to lead to a larger difference, as compared to Burgers' equation. As SGS energy does not dissipate in the KdV equation, but flows back into the resolved scales, the dependence on the true $\mathbf{s}$ is likely higher than for Burgers'.

Overall, the conclusion is that our proposed SP closure model leads to more robust simulations than the CNN, while still improving upon NC with roughly an order of magnitude. 

\subsubsection{Computation time}

Looking at the relative computation times in Figure \ref{fig:convergence} we find that both SP and the CNN are at least $2\times$ faster than the DNS, for Burgers' equation. Furthermore, the computation time for the CNN is slightly shorter, but seems to increase at a larger rate when increasing the $\text{DOF}$.
On our laptop CPU the computation time of a single DNS took no longer than 5 seconds. For such small simulations we expect the computation time to be largely determined by computational overhead. We therefore expect the relative speedup to increase for larger systems.
In fact, the training time amounts to roughly 20 minutes for each neural network. This is rather long compared to a DNS, however for more complex/larger systems we expect the significance of the training time to decrease, relative to the DNS time. Especially for 2D/3D problems where the curse of dimensionality sets in, significant speedups can be expected \cite{List,kochkov}. Also, once a model is trained it can be applied to unseen simulation conditions, without retraining.

For the KdV equation, the relative speedup is roughly $20\times$. This is likely caused by the small time step size required for the DNS ($\Delta t = 10^{-4}$). Larger time step sizes for the DNS consistently resulted in unstable simulations within a few time steps. In addition, we find that computing the true $\mathbf{s}$ for the initial condition does not really affect the computation time.


\subsection{Consistency of the training procedure}
\begin{figure}[ht]
    \centering
    \begin{subfigure}[b]{0.334\textwidth}
        \includegraphics[width = \textwidth]{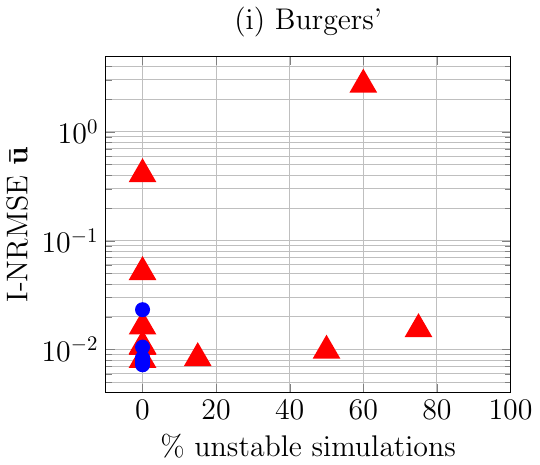}
    \end{subfigure}
    \begin{subfigure}[b]{0.313\textwidth}
        \includegraphics[width = \textwidth]{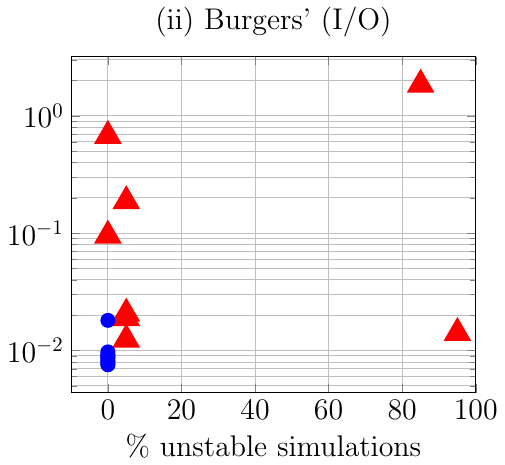}
    \end{subfigure}
    \begin{subfigure}[b]{0.313\textwidth}
        \includegraphics[width = \textwidth]{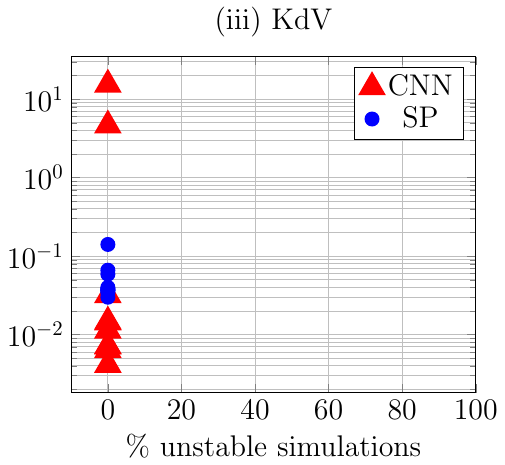}
    \end{subfigure}
    \caption{Integrated solution error evaluated at $T=10$ averaged over 20 simulations and \% of unstable simulations for each closure model in the trained ensemble of 10 closure models ($\text{DOF}=60$). Use cases (i)-(iii) are considered. For (ii) two CNN closure models produced 100\% unstable simulations and are therefore omitted from the graph. $\textbf{CNN}=\text{convolutional neural network closure}$ and $\textbf{SP}=\text{structure-preserving closure}$.}
    \label{fig:ensemble}
\end{figure}

It is important to note that the closure models trained in the previous section possess a degree of randomness. This is caused by the (random) initialization of the network weights and the random selection of the mini-batches. This likely caused the irregular convergence behavior shown in the previous section. In order to evaluate this effect, we train 10 separate replica models for $\text{DOF}=60$, which only differ in the random seed. The trained models are evaluated in terms of stability (number of unstable simulations) and integrated solution error. A simulation is considered unstable when it produces NaN values for $\bar{\mathbf{u}}$. In total 20 simulations per closure model are carried out using the same simulation conditions as in the convergence study. The results are depicted in Figure \ref{fig:ensemble}.
With regard to stability we observe that all trained SP closure models produced exclusively stable simulations. This is in accordance with the earlier derived stability condition \eqref{eq:en_cons_closure}. For the non-periodic test case (ii) we also observe a clear stability advantage, as compared to the CNN.

Regarding this integrated solution error, we observe that the SP closure models all perform very consistently (errors are almost overlapping). The CNNs sometimes outperform SP, but also show very large outliers. This confirms our conclusion of the previous section that our SP closure models are much more robust than the CNNs, which can be `hit or miss' depending on the randomness in the training procedure. However, we still find that SP is often outperformed by the CNN for test case (iii).

\subsection{Structure preservation}

\begin{figure}[ht]
    \centering
    \begin{subfigure}[b]{0.43\textwidth}
        \includegraphics[width = \textwidth]{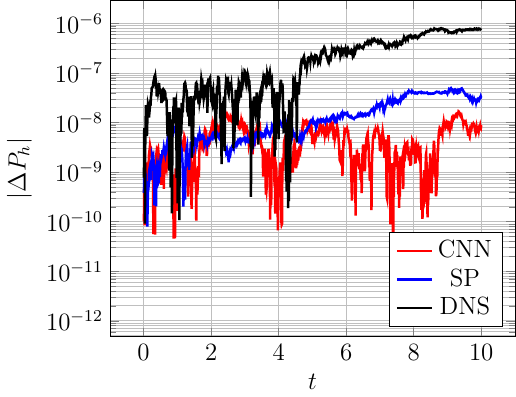}
    \end{subfigure}
    \begin{subfigure}[b]{0.43\textwidth}
        \includegraphics[width = \textwidth]{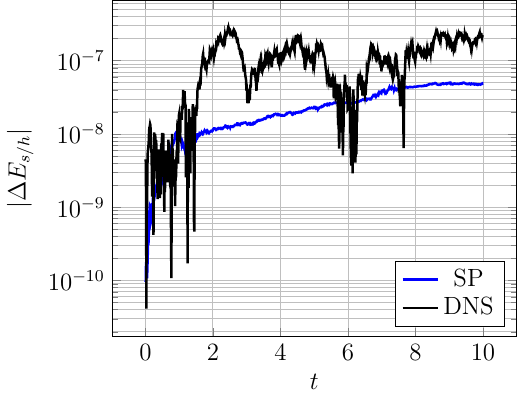}
    \end{subfigure}
    \caption{Change in momentum $\Delta P_h(t) = P_h(t)-P_h(0)$ (left) and total energy $\Delta E_{h/s}(t) = E_{h/s}(t)-E_{h/s}(0)$ (right) for a simulation of the KdV equation with periodic BCs starting from an unseen initial condition. The presented results correspond to $\text{DOF}=40$. $\textbf{DNS}=\text{direct numerical simulation}$, $\textbf{CNN}=\text{convolutional neural network closure}$ and $\textbf{SP}=\text{structure-preserving closure}$.}
    \label{fig:KdV_structure}
\end{figure}

As we formulated a structure-preserving closure model it is important to evaluate if the structure is indeed preserved. For this we consider a single simulation of the KdV equation with periodic BCs, for $\text{DOF}=40$. This is an interesting test case, as the energy should be exactly conserved. In Figure \ref{fig:KdV_structure} we look at both the change in momentum and energy during the simulation.
For momentum conservation we also include the CNN, as this satisfies momentum equation through multiplication by a discrete derivative operator, see \eqref{eq:CNN_closure}. Regarding momentum conservation we find that both SP and the CNN indeed conserve momentum up to machine precision (single-precision). For energy conservation we observe the same for SP. From this we conclude that SP indeed preserves the relevant structure.

\begin{figure}[ht]
\centering
    \begin{subfigure}[b]{0.50\textwidth}
        \includegraphics[width = \textwidth]{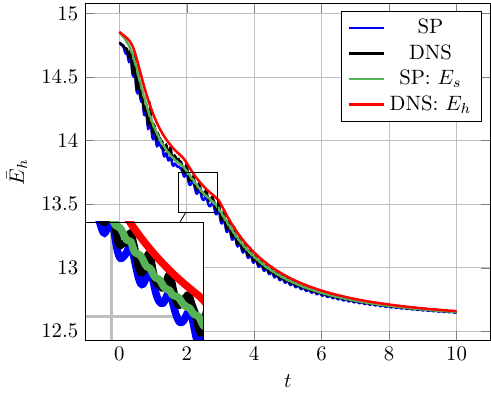}
    \end{subfigure}
    \caption{Resolved and total energy for a simulation of Burgers' equation with periodic BCs starting from an unseen initial condition. The presented results correspond to $\text{DOF}=40$. $\textbf{DNS}=\text{direct numerical simulation}$ and $\textbf{SP}=\text{structure-preserving closure}$.}
    \label{fig:burgers_SP_example}
\end{figure}

Furthermore, we also consider a single simulation of Burgers' equation with periodic BCs, see Figure \ref{fig:burgers_SP_example}.
Here we find that the total energy is always decreasing for SP. However, we find the resolved energy to oscillate due to backscatter, which matches the filtered DNS result.
From this we conclude that, although the energies do not match the DNS result exactly, SP indeed allows for backscatter to be modeled correctly.

\subsection{Burgers' equation \& energy spectra}

\begin{figure}[ht]
\begin{subfigure}[b]{0.3332\textwidth}
        \includegraphics[width = \textwidth]{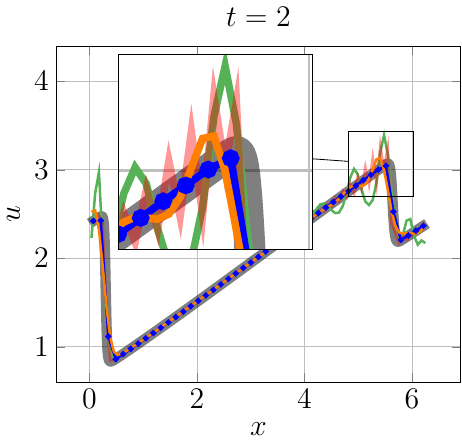}
    \end{subfigure}
    \begin{subfigure}[b]{0.3134\textwidth}
        \includegraphics[width = \textwidth]{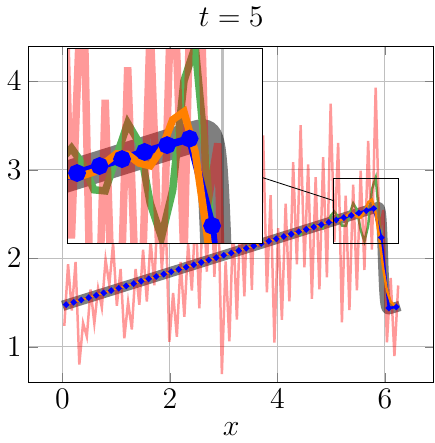}
    \end{subfigure}
    \centering
    \begin{subfigure}[b]{0.3134\textwidth}
        \includegraphics[width = \textwidth]{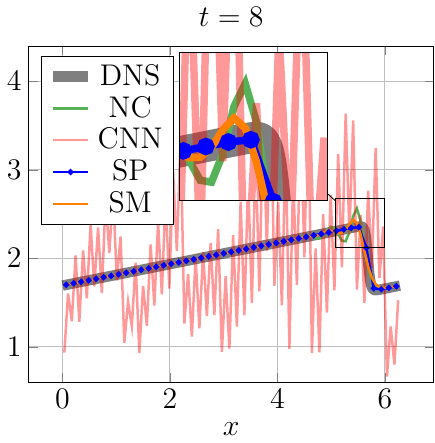}
    \end{subfigure}
    \caption{Solutions of Burgers' equation with periodic BCs at different points in time, starting from an unseen initial condition. The solutions are produced by the different closure models corresponding to $\text{DOF}=90$. $\textbf{DNS}=\text{direct numerical simulation}$, $\textbf{NC}=\text{no closure}$, $\textbf{CNN}=\text{convolutional neural network closure}$, $\textbf{SP}=\text{structure-preserving closure}$, and $\textbf{SM}=\text{constant Smagorinsky model}$.}
    \label{fig:burgers}
\end{figure}

To find out what happens when a CNN becomes unstable we consider Burgers' equation with periodic BCs, for $\text{DOF}=90$. For this CNN we observed a rather large error spike in the convergence study (Figure \ref{fig:convergence}). To analyze this error spike we consider a single simulation starting from an unseen initial condition. The results are depicted in Figure \ref{fig:burgers}. Here we observe a large buildup of numerical noise for the CNN, as the simulation progresses. In addition, we observe that SP nicely suppresses the wiggles produced by NC around the shock \cite{Jameson}. SM also manages to do this, although to a lesser extent.

\begin{figure}[ht]
\begin{subfigure}[b]{0.43\textwidth}
        \includegraphics[width = \textwidth]{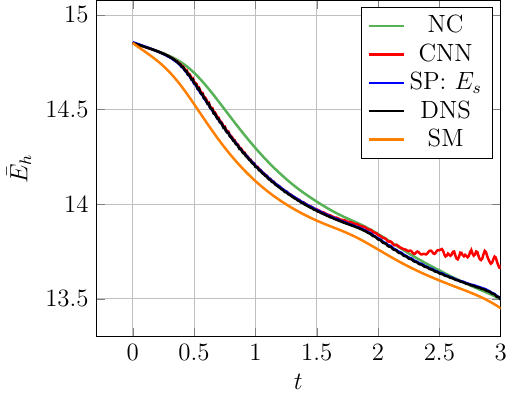}
    \end{subfigure}
    \begin{subfigure}[b]{0.388\textwidth}
        \includegraphics[width = \textwidth]{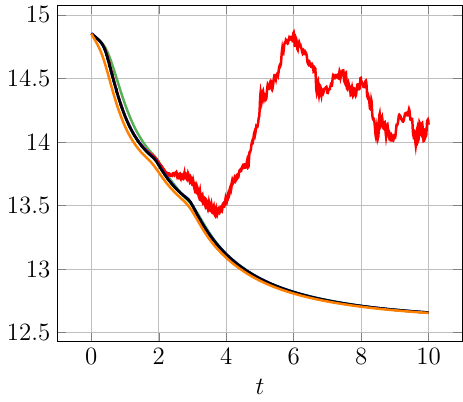}
    \end{subfigure}
    \centering
    \begin{subfigure}[b]{0.43\textwidth}
        \includegraphics[width = \textwidth]{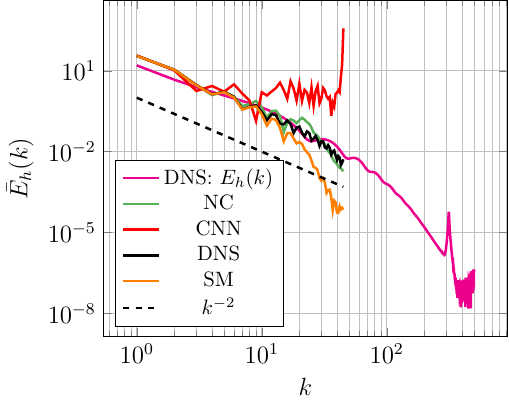}
    \end{subfigure}
    \begin{subfigure}[b]{0.396\textwidth}
        \includegraphics[width = \textwidth]{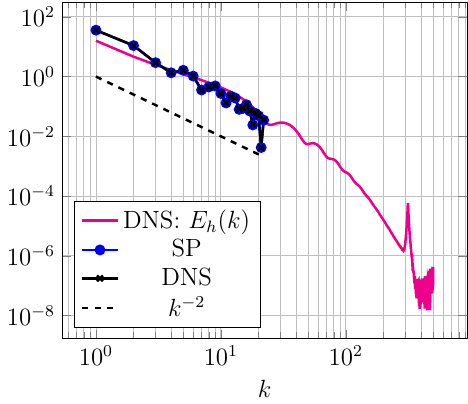}
    \end{subfigure}
    
    \caption{(Top) Resolved energies corresponding to the simulation in Figure \ref{fig:burgers} using the different closures for $\text{DOF}=90$. For SP we depicted the total energy. The SP energy trajectory mostly overlaps with the DNS. (Bottom) The average energy spectra of this simulation, evaluated on the interval $3\leq t\leq 7$. The left plot corresponds to $I = \text{DOF}$ and the right plot to $I=\text{DOF}/2$. Furthermore, we also show the spectra of the unfiltered DNS. These slightly differ from the filtered DNS in the small wavenumbers, as our filter is not a spectral filter \cite{Edeling}. $\textbf{DNS}=\text{direct numerical simulation}$, $\textbf{NC}=\text{no closure}$, $\textbf{CNN}=\text{convolutional neural network closure}$, $\textbf{SP}=\text{structure-preserving closure}$, and $\textbf{SM}=\text{constant Smagorinsky model}$.}
    \label{fig:burgers_energy}
\end{figure}

Next, we look at the energy trajectories and energy spectra corresponding to this simulation, see Figure \ref{fig:burgers_energy}. The energy spectra are given as a function of the wavenumber $k$. These spectra are computed by carrying out a discrete Fourier transform of $\bar{\mathbf{u}}$ and computing the energy for every $k$ \cite{Edeling}. Furthermore, the energy spectra are only depicted up to the wavenumber resolved by the closure models.
We find that for the CNN the energy starts to diverge at around $t = 2$. This worsens as the simulation progresses and a large increase in energy is observed. Looking at the energy spectrum this corresponds to a buildup of energy in the small scales (large $k$) and the numerical noise observed in the solution. This can cause the simulation to blow up.  Furthermore, we find that SM is too dissipative. This results mostly in a lack of energy in the small scales. In addition, NC is not dissipative enough, resulting in a slight buildup of energy in the small scales. Finally, SP seems to capture the energy balance nicely, looking at both the trajectory and the spectrum. SP, and to a lesser extent NC, show the expected $k^{-2}$ slope in the energy spectrum \cite{k2_slope}.

\subsection{Extrapolation in parameter space}

\begin{figure}[ht]
    \centering
        \includegraphics[width = 0.50\textwidth]{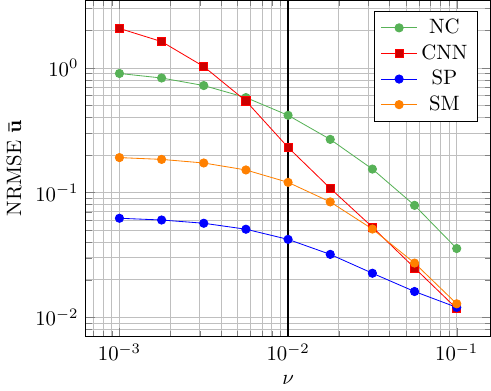}
    \caption{Integrated solution error evaluated at $T=5$ averaged over the three trained closure models and three simulations for each viscosity value. The black vertical separates the training range (right) from the extrapolation range (left). $\textbf{NC}=\text{no closure}$, $\textbf{CNN}=\text{convolutional neural network closure}$, $\textbf{SP}=\text{structure-preserving closure}$, and $\textbf{SM}=\text{constant Smagorinsky model}$.}
    \label{fig:visc_extrapolation}
\end{figure}

Next, we are interested in how well such closure models are capable of extrapolating in parameter space. In particular, we consider different viscosity values $\nu$ for the Burgers' equation. For this purpose we provide $\nu$ as an additional input to the neural network. For the training data we consider $\nu \in [10^{-2},10^{-1.75},10^{-1.5},10^{-1.25},10^{-1}]$. To generate the training data we randomly generate three different initial conditions and carry out a DNS for each considered $\nu$. To accommodate smaller values of $\nu$ used in the extrapolation experiment (up to $\nu = 10^{-3}$) we refine the grid to $N=5000$. Furthermore, we use $\Delta t = 10^{-5}$ for the time integration and simulate up to $T = 5$. The smaller time step size is used to accommodate for the finer grid. Every 1000th time step we save the solution and use this as training data for the closure models. For both SP and the CNN we train three closure models which only differ in the random seed used for the training procedure and initial parameter values. The networks underlying the methodologies are kept small, namely only a single hidden layer with 20 channels. This is done, as smaller networks show higher potential for generalization \cite{maulikshankar2022differentiable}. The remaining hyperparameters and the training procedure are kept the same, see \ref{sec:training}. For our closure models we reduce the $\text{DOF}$ of the system from $5000$ to $50$ (a 100 fold reduction). We also use a larger time step size of $\overline{\Delta t}=0.01$ (a 1000 fold increase). We evaluate the performance of each of the closure models for three simulations starting from different initial conditions. For each initial condition we consider 9 different values of $\nu$ in the range $[10^{-3},10^{-1}]$. The integrated solution error, see \eqref{eq:int_sol}, is once again used as the performance metric. The results are depicted in Figure \ref{fig:visc_extrapolation}.
For each model we find that the error increases for smaller values of $\nu$. This is to be expected as smaller values of $\nu$ lead to stronger spatial gradients (shocks) which are harder to deal with in a numerical setting \cite{Jameson,benjamin_thesis}. When comparing the performance of the different closures we find that our proposed SP scheme consistently outperforms the other closure models. This is the case for both inter- and extrapolation. Regarding the CNN we observe a drop in performance as the viscosity decreases, even within the training range. It seems that the CNN is not able to adapt to the range of values of $\nu$, as previously the CNN was able to perform quite well for $\text{DOF}=50$ and a single viscosity value, see Figure \ref{fig:convergence}. 
Furthermore, we find that SM outperforms NC across the considered range. This means that the obtained Smagorinsky constant generalizes well. This is to be expected as it is a simple model containing only a single parameter \cite{smagorinsky1963general,maulikshankar2022differentiable}.

\subsection{Extrapolation in space and time}\label{sec:extrapolation}

As a final test case we evaluate how well the closure models are capable of extrapolating in space and time. We consider the KdV equation on an extended spatial domain $\Omega = [0,96]$ and run the simulation until $T = 50$. This corresponds to a $3\times$ and $5\times$ increase, respectively, with respect to the training data. We choose the KdV equation for this experiment as it is non-dissipative. It therefore leads to more interesting long-term behavior. As closure models, we use the ones trained during the convergence study that correspond to the grid-spacing of the grids employed in this test case. 

\begin{figure}[ht]
    \centering
    \begin{subfigure}[b]{0.332\textwidth}
        \includegraphics[width = \textwidth]{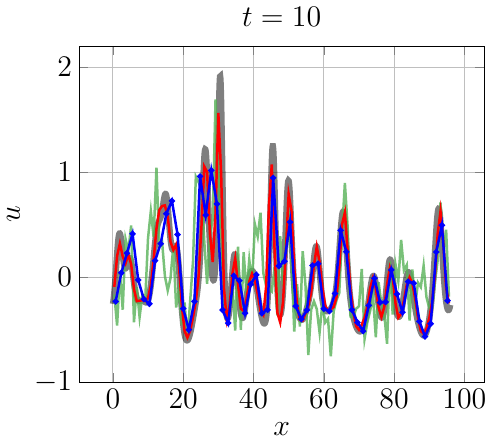}
    \end{subfigure}
    \begin{subfigure}[b]{0.314\textwidth}
        \includegraphics[width = \textwidth]{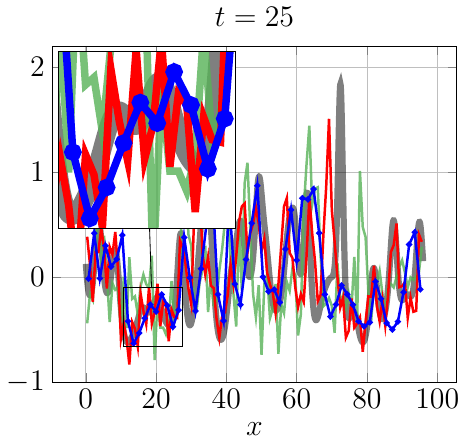}
    \end{subfigure}
    \begin{subfigure}[b]{0.314\textwidth}
        \includegraphics[width = \textwidth]{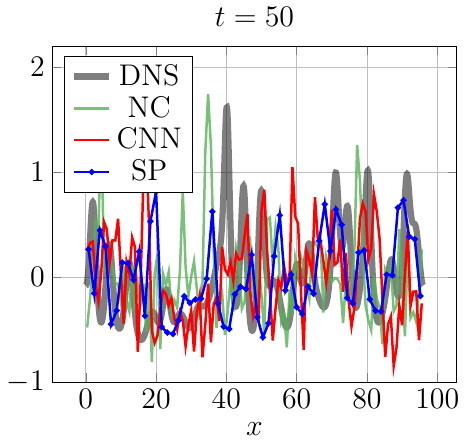}
    \end{subfigure}
    \caption{Solutions of the KdV equation with periodic BCs at different point in time, starting from an unseen initial condition. The solutions are produced by the different closure models corresponding to $\text{DOF}=40$.  In this case, the spatial and temporal domain are increased $3\times$ and $5\times$ with respect to the training data. $\textbf{DNS}=\text{direct numerical simulation}$, $\textbf{NC}=\text{no closure}$, $\textbf{CNN}=\text{convolutional neural network closure}$, and $\textbf{SP}=\text{structure-preserving closure}$.}
    \label{fig:KdV_large}
\end{figure}

Let us start by looking at snapshots from different points in the simulation, see Figure \ref{fig:KdV_large}. At the edge of the training region $t=10$ we find that the solutions are still roughly aligned, except for NC which already contains a lot of numerical noise. At $t=25$ we find that both SP and the CNN start to diverge from the DNS. In addition, we observe a buildup of numerical noise for the CNN. This worsens at $t=50$, while SP seems to remain free of numerical noise.

\begin{figure}[ht]
    \centering
        \includegraphics[width = 0.50\textwidth]{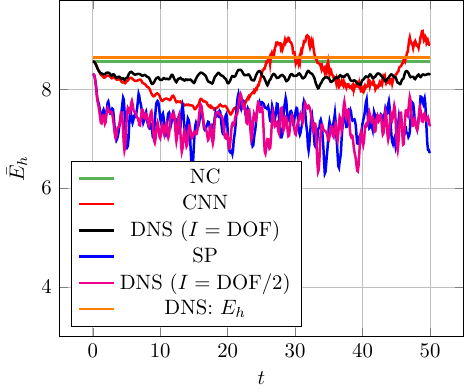}
    \caption{Trajectory of the resolved energy $\bar{E}_h$ for the simulation presented in Figure \ref{fig:KdV_large} for each of the different models corresponding to $\text{DOF}=40$. The DNS resolved energy is depicted for both $I=\text{DOF}$ (to compare with the CNN) and $I = \text{DOF}/2$ (to compare with SP). $\textbf{DNS}=\text{direct numerical simulation}$, $\textbf{NC}=\text{no closure}$, $\textbf{CNN}=\text{convolutional neural network closure}$, and $\textbf{SP}=\text{structure-preserving closure}$.}
    \label{fig:extended_E_bar}
\end{figure}

To make a more thorough analysis we consider the trajectories of the resolved energy. This is presented in Figure \ref{fig:extended_E_bar}.
We find that for SP the resolved energy (in blue) stays in close proximity to the filtered DNS (in magenta). This is in contrast to the CNN (in red) which starts to diverge from the DNS (in black) around $t=5$. The resolved energy for the CNN also exceeds the maximum allowed total energy $E_h$ (in orange) at different points in the simulation, which is nonphysical. We conclude that adding the SGS variables and conserving the total energy helps with capturing the delicate energy balance between resolved and SGS energy that characterizes the DNS. 
It is also interesting to note that NC (in green) conserves the resolved energy, as the coarse discretization conserves the discrete energy. 

\begin{figure}[ht]
    \centering
    \begin{subfigure}[b]{0.334\textwidth}
        \includegraphics[width = \textwidth]{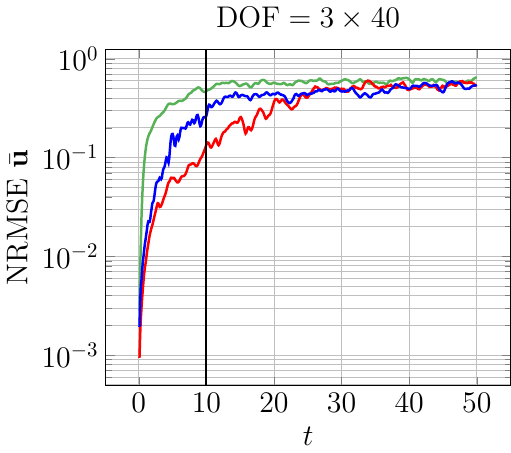}
    \end{subfigure}
    \centering
    \begin{subfigure}[b]{0.313\textwidth}
        \includegraphics[width = \textwidth]{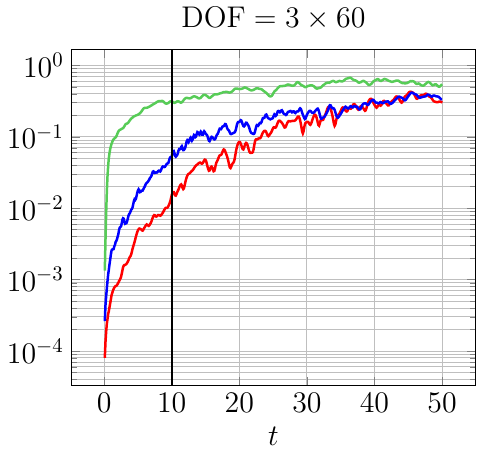}
    \end{subfigure}
    \centering
    \begin{subfigure}[b]{0.313\textwidth}
        \includegraphics[width = \textwidth]{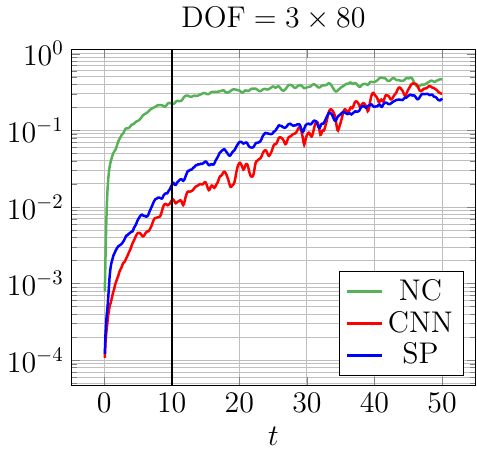}
    \end{subfigure}
    \centering
    \begin{subfigure}[b]{0.328\textwidth}
        \includegraphics[width = \textwidth]{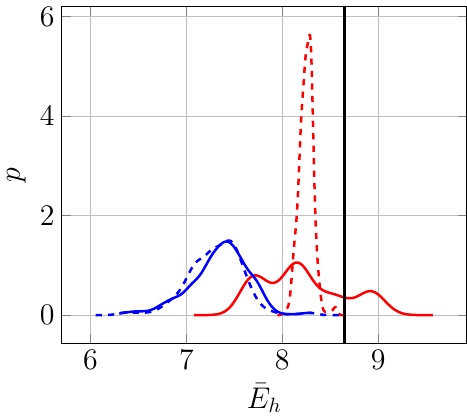}
    \end{subfigure}
    \centering
    \begin{subfigure}[b]{0.316\textwidth}
        \includegraphics[width = \textwidth]{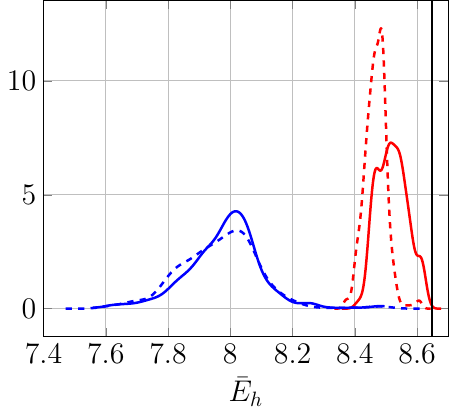}
    \end{subfigure}
    \centering
    \begin{subfigure}[b]{0.316\textwidth}
        \includegraphics[width = \textwidth]{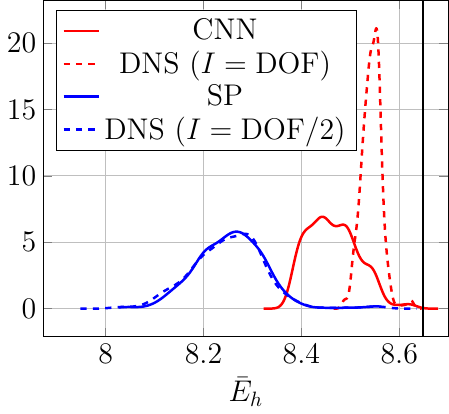}
    \end{subfigure}
    \caption{Solution error trajectory (top) and KDEs estimating the distribution of $\bar{E}_h$ (bottom) for the trained closure models corresponding to different numbers of $\text{DOF}$. These quantities are computed for a simulation of the KdV equation with the same initial condition on the extended spatial and temporal domain. In the top row the vertical black line indicates the maximum time present in the training data, while in the bottom row it indicates the total energy of the DNS (which should not be exceeded). The DNS resolved energy is again depicted for both $I=\text{DOF}$ (to compare with the CNN) and $I = \text{DOF}/2$ (to compare with SP). $\textbf{NC}=\text{no closure}$, $\textbf{CNN}=\text{convolutional neural network closure}$, and $\textbf{SP}=\text{structure-preserving closure}$.}
    \label{fig:err_KDE}
\end{figure}

To make a more quantitative analysis of this phenomenon we investigate the trajectory of the solution error and the Gaussian kernel density estimate (KDE) \cite{KDE} of the resolved energy distributions. The latter analysis is carried out to analyze whether the closure models capture the correct energy balance between the resolved and SGS energy. The results for $\text{DOF} \in \{40,60,80\}$ are depicted in Figure \ref{fig:err_KDE}.
Looking at the solution error trajectories we find that at the earlier stages of the simulation the CNN outperforms SP.
However, SP slowly catches up with the CNN past the training region. 
With regards to the resolved energy distribution we find that for each of the considered numbers of $\text{DOF}$ SP is capable reproducing the DNS distribution. On the other had, the CNN closure models struggle to capture this distribution. For $\text{DOF}=40$ a significant part of the distribution even exceeds the total energy present in the DNS, i.e. there occurs a nonphysical influx of energy. 

From this we conclude that both SP and the CNN closure models are capable of extrapolating beyond the training data. However, only SP is capable of correctly capturing the energy balance between the resolved and unresolved scales. This allows it to more accurately capture the statistics of the DNS result.

\section{Discussion on the applicability to 2D/3D Navier-Stokes}\label{sec:discussionNS}

A number of challenges arises when considering to apply the proposed methodology to the 2D/3D Navier-Stokes equations. Firstly, the discrete nature of the presented closure modeling framework, making it dependent on the compression factor $J$ between the coarse and the fine-grid resolution. This is a result of taking the `discretize first' approach. While this allows the method to be highly specialized to the discretization, it can limit the applicability of the closure model to grid pairs with the associated $J$. One possibility is to train a single closure model with training data from multiple compression factors. Another possibility is to consider a continuous formulation of the closure problem, along with a continuous expression for the SGS variables. The finite element framework might be a useful starting point, as it approximates the solution using a continuous, but still highly local, basis \cite{girault2012finite}. In addition, it easily allows for non-Cartesian grids. The convolutional layers in the closure model could then be replaced by graph convolutions \cite{graph_NN_belbute2020combining}, which work for unstructured grids. Another useful alternative is to consider the Fourier neural operator, which promises to be discretization invariant \cite{fourier_neural_operators}. Furthermore, one could do a sparse regression on the closure model to identify physical terms \cite{BRUNTON2016710_sindy}, and train the model for different $J$ to understand the effect of discretization error. This would effectively separate the two sources of error.
Note that our framework can also be applied to nonuniform Cartesian grids, as long as the compression factor $J$ is constant, and an energy-conserving discretization is available. The size of the computational domain can even be increased by an arbitrary factor as compared to the training domain, without retraining the closure model. The latter was shown in section \ref{sec:extrapolation}. For even more irregular grids graph neural networks could be used to increase the flexibility of the framework \cite{graph_NN_belbute2020combining}.

A second challenge is the linear compression for $\mathbf{s}$, which might be a limiting factor for 2D/3D flows. As an alternative, non-linear compression methods like autoencoders could be employed \cite{bank2023autoencoders}. Furthermore, the natural extension of the compressed SGS energy $\mathbf{s}$ to 2D and 3D is probably the sub-grid scale stress tensor. This means that our framework needs to be extended with multiple sub-grid scale equations, e.g.\ $\mathbf{s}_{1}$, $\mathbf{s}_{2}$, $\ldots$
Finally, we could draw inspiration from the scale-adaptive simulation method which also uses the SGS energy in its equations \cite{menter2010_SAS_scale}.

\section{Conclusion}\label{sec:conclusion}

In this paper we proposed a novel way of constructing machine learning-based closure models in a structure-preserving fashion.
We started by applying a spatial averaging filter to a fine-grid solution and writing the resulting system in closure model form. We showed that by applying this filter we effectively remove part of the energy. Next, we introduced a linear compression of the subgrid-scales (SGSs) into a set of SGS variables, defined on the coarse grid. These serve as a means of reintroducing the removed energy back into the system. This allows us to use the concept of kinetic energy conservation in closure modeling. In turn we introduced an extended system of equations which models the evolution of the filtered solution as well as the evolution of the SGS variables. For this extended system we proposed a structure-preserving closure modeling framework which allows for energy exchange between the filtered solution and the SGS variables, in addition to dissipation. This framework serves to constrain the underlying convolutional neural network (CNN) such that no additional energy enters the system. In this way we achieve stability by abiding by the underlying energy conservation law. The advantage is that the framework still allows for backscatter through the energy present in the SGS variables. In addition, momentum conservation is also satisfied. Finally, the framework was applied to both Burgers' and Korteweg-de Vries (KdV) equation.

A convergence study showed that the learned SGS variables are able to accurately match the original SGS energy content, with accuracy consistently improving when refining the coarse-grid resolution. 

Given the SGS compression, our proposed structure-preserving framework (SP) was compared to a vanilla CNN. Overall, our SP method performed roughly on par with the CNN in terms of accuracy, albeit the CNN outperformed SP slightly for the KdV equation. However, the results for the CNN were typically inconsistent, not showing clear convergence of the error upon increasing the grid resolution. In addition, \textit{the CNN suffered from stability issues}. On the other hand, our SP method produced stable results in all cases, while also consistently improving upon the `no closure model' result. To be more specific, it did so by roughly an order of magnitude in terms of reproducing the reference solution. 

This conclusion was further strengthened by training an ensemble of closure models. This was done to investigate the consistency of the closure model performance with respect to the randomness inherent to the training procedure. We observed that the trained vanilla CNNs differed significantly in performance and stability, whereas the different SP models performed very similarly to each other. The SP closures also displayed no stability issues. Our SP framework has therefore shown to be more robust and successfully resolves the stability issues which plague conventional CNNs. 

Our numerical experiments confirmed the structure-preserving properties of our method: exact momentum conservation, energy conservation up to a time discretization error, and a strictly decreasing energy in the presence of dissipation. We also showed that our method succeeds in accurately modeling backscatter in both Burgers' and the KdV equation. 

Regarding extrapolation in parameter space, out method outperformed both the CNN and the Smagorinsky model on viscosity values in- and outside the training range. From this we conclude our methodology leads to more accurate results than the conventional methods for this type of application.
Furthermore, when extrapolating in space and time, the advantage of including the SGS variables and embedding structure-preserving properties became even more apparent: Our method is much better at capturing the delicate energy balance between the resolved and SGS energy. This in turn yielded better long-term error behavior. 

Based on these results we conclude that including the SGS variables, as well as adherence to the physical structure, has the important advantages of stability and long-term accuracy. In addition, it also leads to more consistent performance. This work therefore serves as an important starting point for building physical constraints into machine learning-based turbulence closure models. More generally, our framework is potentially applicable to a wide range of systems that possess multiscale behavior, while also possessing a secondary conservation law, for example incompressible pipe flow \cite{buist2021energy}. Currently our efforts are mainly directed towards the incompressible Navier-Stokes equations. For instance, Kolmogorov flow would be a good starting point \cite{chandler_kerswell_2013_kolmogorov}.

\section*{CRediT authorship contribution}

\textbf{T. van Gastelen:} Conceptualization, Methodology, Software, Writing - original draft. \textbf{W. Edeling:} Writing - review \& editing. \textbf{B. Sanderse:} Conceptualization, Methodology, Writing - review \& editing, Funding acquisition.

\section*{Data availability}

The code used to generate the training data and the implementation of the neural networks can be found at \url{https://github.com/tobyvg/ECNCM_1D}.

\section*{Acknowledgements}

This publication is part of the project ``Unraveling Neural Networks with Structure-Preserving Computing” (with project number OCENW.GROOT.2019.044
of the research programme NWO XL which is financed by the Dutch Research
Council (NWO)). Part of this publication is funded by Eindhoven University of Technology. Finally, we thank the reviewers for their feedback, enhancing the quality of the article.
\clearpage

\bibliographystyle{elsarticle-num}

\bibliography{references} 
\appendix

\clearpage

\section{Filter properties}\label{sec:filter_properties}

Here we derive important properties of the spatial averaging filter. We first show that equation \eqref{eq:filtering_integ_prop} holds:
\begin{equation}
    (\mathbf{R}\bar{\mathbf{a}},\mathbf{R}\bar{\mathbf{b}})_\omega = \bar{\mathbf{a}}^T\mathbf{R}^T\boldsymbol{\omega}\mathbf{R}\bar{\mathbf{b}} = \bar{\mathbf{a}}^T\boldsymbol{\Omega} \mathbf{W}\boldsymbol{\omega}^{-1}\boldsymbol{\omega}\mathbf{R}\bar{\mathbf{b}}= \bar{\mathbf{a}}^T\boldsymbol{\Omega} \underbrace{\mathbf{W}\mathbf{R}}_{=\mathbf{I}}\bar{\mathbf{b}} = (\bar{\mathbf{a}},\bar{\mathbf{b}})_\Omega,
\end{equation}
where we used that fact that
\begin{equation}
    \boldsymbol{\Omega}\mathbf{W}\boldsymbol{\omega}^{-1} = \mathbf{R}^T.
\end{equation}
Next, we proof that $\mathbf{R}\bar{\mathbf{u}}$ is orthogonal to $\mathbf{u}^\prime$:
\begin{equation}
    \begin{split}
        (\mathbf{R}\bar{\mathbf{u}},\mathbf{u}^\prime)_\omega &= (\mathbf{R}\bar{\mathbf{u}},\mathbf{u}-\mathbf{R}\bar{\mathbf{u}})_\omega = (\mathbf{R}\bar{\mathbf{u}},(\mathbf{I}-\mathbf{R}\mathbf{W})\mathbf{u})_\omega = \bar{\mathbf{u}}^T\mathbf{R}^T\boldsymbol{\omega}(\mathbf{I}-\mathbf{R}\mathbf{W})\mathbf{u}   \\&= \bar{\mathbf{u}}^T \boldsymbol{\Omega}\mathbf{W}(\mathbf{I}-\mathbf{R}\mathbf{W})\mathbf{u} = (\bar{\mathbf{u}}, (\mathbf{W}-\underbrace{\mathbf{W}\mathbf{R}}_{=\mathbf{I}}\mathbf{W})\mathbf{u})_\Omega = (\bar{\mathbf{u}} ,(\mathbf{W}-\mathbf{W})\mathbf{u})_\Omega = 0.
    \end{split}
\end{equation}
Finally, we show that equation \eqref{eq:mom_filter} holds:
\begin{equation}
    (\mathbf{1}_\omega,\mathbf{u})_\omega =  \mathbf{1}_\omega^T\boldsymbol{\omega}\mathbf{u}  =\mathbf{1}_{\Omega}^T\mathbf{R}^T\boldsymbol{\omega}\mathbf{u}= \mathbf{1}_\Omega^T\boldsymbol{\Omega}\mathbf{W}\boldsymbol{\omega}^{-1}\boldsymbol{\omega}\mathbf{u} =\mathbf{1}^T_\Omega\boldsymbol{\Omega}\mathbf{W}\mathbf{u}  =(\mathbf{1}_\Omega,\mathbf{W}\mathbf{u})_\Omega = (\mathbf{1}_\Omega,\bar{\mathbf{u}})_\Omega,
\end{equation}
where we used the fact that $\mathbf{1}_\omega = \mathbf{R}\mathbf{1}_\Omega$.
\section{Comparing coarse and fine-grid dissipation}\label{sec:diss-diff}
\begin{figure}[h]
    \centering
     \includegraphics[width = 0.55\textwidth]{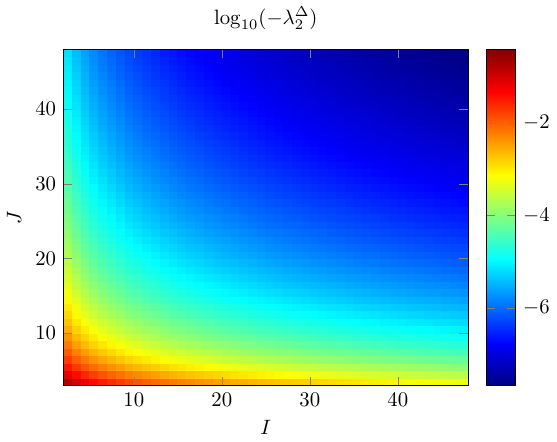}
    \caption{Largest non-zero eigenvalue $\lambda^\Delta_2$ of $\mathbf{D}_\Delta$ for different values of $I$ and $J$.}
    \label{fig:eigenvalues_D2}
\end{figure}

Here we compare the rate of dissipation induced by the 1D diffusion operator discretized on the coarse grid, $\bar{\mathbf{D}} \in \mathbb{R}^{I\times I}$ (for grid-spacing $H$), with the dissipation induced by the same operator but discretized on the fine grid, $\mathbf{D} \in \mathbb{R}^{N\times N}$ (for grid-spacing $h=\frac{H}{J}$). The difference in dissipated energy between these two quantities $\Delta_D$ is given by:
\begin{equation}\label{eq:diss-diff}
\begin{split}
    \Delta_D &= (\mathbf{u},\mathbf{D} \mathbf{u})_\omega - (\bar{\mathbf{u}},\bar{\mathbf{D}}\bar{\mathbf{u}})_\Omega = \mathbf{u}^T\boldsymbol{\omega} \mathbf{D}\mathbf{u} - (\mathbf{W}\mathbf{u})^T\boldsymbol{\Omega}\bar{\mathbf{D}}\mathbf{W}\mathbf{u} =  \mathbf{u}^T \frac{H}{J}\mathbf{D}\mathbf{u} -   \mathbf{u}^TH\mathbf{W}^T \bar{\mathbf{D}}\mathbf{W}\mathbf{u} \\
    &= \frac{J}{H} \mathbf{u}^T  (\frac{H^2}{J^2}\mathbf{D} - \frac{H^2}{J} \mathbf{W}^T \bar{\mathbf{D}}\mathbf{W})\mathbf{u} = \frac{J}{H} \mathbf{u}^T \underbrace{ (h^2\mathbf{D} - \frac{1}{J} \mathbf{W}^T (H^2 \bar{\mathbf{D}})\mathbf{W})}_{=:\mathbf{D}_\Delta}\mathbf{u},
\end{split}
\end{equation}
for periodic BCs.
Here we have written $\mathbf{D}^\Delta$ such that it is independent of the grid-spacing, but only depends on the ratio $J = \frac{H}{h}$. Note that $\mathbf{D}^\Delta$ is a symmetric matrix and as such its eigenvalues $\lambda_N^\Delta \leq \ldots \leq \lambda_1^\Delta$ are real and its eigenvectors can be chosen to form an orthogonal basis. Furthermore, $\lambda_{1}^\Delta=0$  for periodic boundary conditions. We can bound $\Delta_D$ by noting that
\begin{equation}
    \Delta_D \leq \max_i \frac{J}{H}\lambda^\Delta_i ||\mathbf{u}||_2^2.
\end{equation}
If the eigenvalues are all strictly nonpositive, it follows that the difference in dissipation $\Delta_D$ is always less than or equal to zero. In other words, $\mathbf{D}$ extracts more (or equal) energy from the reference system as $\bar{\mathbf{D}}$ does from the filtered system. To prove that the eigenvalues of $\mathbf{D}_\Delta$ are indeed nonpositive turns out to be a difficult problem which we circumvent with a numerical `proof'. In Figure \ref{fig:eigenvalues_D2} we display the largest non-zero eigenvalue $\lambda^\Delta_2$ for different values of $I$ and $J$, indicating that $\lambda^{\Delta}_{i}\leq 0$ for realistic values of $I$ and $J$.
\section{SGS compression}\label{sec:SGS_compression}

In this section we outline how we obtain the SGS compression parameter values $\mathbf{t} \in \mathbb{R}^J$ such that $\mathbf{s}\approx \mathbf{W}(\mathbf{u}^\prime)^2$. This can be achieved by using a singular value decomposition (SVD). The SVD minimizes the following loss function for $\hat{\mathbf{t}} \in \mathbb{R}^J$ (here we assume a uniform grid):
\begin{equation}\label{eq:L_s}
    \begin{split}
    \mathcal{L}_s(\mathbf{X}_\mu;\boldsymbol{\Theta}) &= \frac{1}{pI}\sum_{\boldsymbol{\mu}\in \mathbf{X}_\mu}|J^{-1} \boldsymbol{\mu}^T \boldsymbol{\mu}-\boldsymbol{\mu}^T\mathbf{t}\mathbf{t}^T\boldsymbol{\mu}| =  \frac{1}{pIJ}\sum_{\boldsymbol{\mu}\in \mathbf{X}_\mu}| \boldsymbol{\mu}^T \boldsymbol{\mu}-\boldsymbol{\mu}^T\hat{\mathbf{t}}\hat{\mathbf{t}}^T\boldsymbol{\mu}|  \\ &= \frac{1}{pIJ}\sum_{\boldsymbol{\mu}\in \mathbf{X}_\mu}|\boldsymbol{\mu}^T\boldsymbol{\mu} + \boldsymbol{\mu}^T\hat{\mathbf{t}}\hat{\mathbf{t}}^T\hat{\mathbf{t}}\hat{\mathbf{t}}^T\boldsymbol{\mu} -  2\boldsymbol{\mu}^T\hat{\mathbf{t}}\hat{\mathbf{t}}^T\boldsymbol{\mu}|= \frac{1}{pIJ}\sum_{\boldsymbol{\mu}\in \mathbf{X}_\mu}(\hat{\mathbf{t}}\hat{\mathbf{t}}^T\boldsymbol{\mu} - \boldsymbol{\mu})^T(\hat{\mathbf{t}}\hat{\mathbf{t}}^T\boldsymbol{\mu} - \boldsymbol{\mu}) \\  \quad &\text{subject to} \quad \hat{\mathbf{t}}^T \hat{\mathbf{t}}= 1,
\end{split}
\end{equation}
where $\boldsymbol{\mu}$ refers to a column vector of snapshot matrix $\mathbf{X}_\mu \in \mathbb{R}^{J \times Ip}$ and $\mathbf{t}= J^{-1/2}\hat{\mathbf{t}}$.
The last expression for $\mathcal{L}_s$ is a projection error. Such an error is typically minimized in a reduced order modeling setting to obtain a proper basis, see \cite{podbenjamin}. For this one typically uses an SVD. Conveniently, minimizing this error is equivalent to minimizing our required energy error (first expression). In this expression the prefactor $J^{-1}$ accounts for $\mathbf{W}$ in the SGS energy \eqref{eq:SGS_energy}. The snapshot matrix $\mathbf{X}_\mu$ is constructed from the training data set $\mathbf{X}_\text{u} \in \mathbb{R}^{N \times p}$ which contains $p$ DNS snapshots as the columns. From this matrix we compute the SGS content: $\mathbf{X}_{\text{u}^\prime} = (\mathbf{I}-\mathbf{R}\mathbf{W})\mathbf{X}_\text{u}$. Finally, this matrix is reshaped into $\mathbf{X}_\mu$ such that each column corresponds to the SGS content in a coarse cell, i.e. $\mathbf{u}^\prime \rightarrow [\boldsymbol{\mu}_i,\ldots,\boldsymbol{\mu}_I]$ for each column in $\mathbf{X}_{\text{u}^\prime}$. The SVD of $\mathbf{X}_\mu$ is given by
\begin{equation}
    \mathbf{X}_\mu = \mathbf{U}_\mu\boldsymbol{\Sigma}_\mu\mathbf{V}_\mu^T,
\end{equation}
where $\mathbf{U}_\mu \in \mathbb{R}^{J \times J}$ and $\mathbf{V}_\mu \in \mathbb{R}^{Ip \times Ip}$ are unitary matrices containing the left and right-singular vectors, respectively, and $\boldsymbol{\Sigma}_\mu \in \mathbb{R}^{J \times Ip}$ contains the singular values on the diagonal. The values for $\hat{\mathbf{t}}$ that minimize \eqref{eq:L_s} correspond to the first column of $\mathbf{U}_\mu$.

\section{Non-periodic boundary conditions}\label{sec:BCs}

To extend our method to different types of BCs we resort to what the machine learning community refers to as padding \cite{conv_NN} and the scientific computing community refers to as the ghost-cell method \cite{ghost_cell_DADONE20071513}. We will treat both inflow and outflow BCs, on uniform 1D grids, as this is relevant for Burgers' equation. 

\subsection{Implementation for the fine grid}

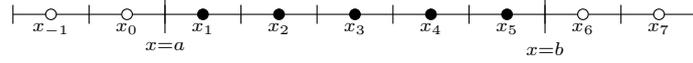
\begin{figure}[ht]
    \centering
    \begin{tikzpicture}
    \draw  (1,-0.125) -- (1,0.125);
    \draw  (3,-0.25) -- (3,0.25);
    \fill[black]   (3,-0.25) circle (0.0cm) node[below]{$\mathsmaller{x=a}$};
    \draw  (8,-0.25) -- (8,0.25);
    \fill[black]   (8,-0.25) circle (0.0cm) node[below]{$\mathsmaller{x=b}$};
    \foreach \i in {1,2,3}{

        \foreach \j in {1,2,3}{
        \def\ij{\the\numexpr 3 * \i - 3 + \j-2}
        \draw  (1+3*\i-3+\j,-0.125) -- (1+3*\i-3+\j,0.125);
        \draw (3*\i-3+\j,0) -- (1+3*\i-3+\j,0);
        \ifthenelse{\ij < 1 \OR \ij > 5}{\draw[black,fill = white]  (1+3*\i-3+\j -0.5,0) circle (0.07cm)  node[below]{$\mathsmaller{ {x_{\ij}}}$};}{\draw[black,fill = black]  (1+3*\i-3+\j -0.5,0) circle (0.07cm)  node[below]{$\mathsmaller{ {x_{\ij}}}$};}
    
        }}
    \end{tikzpicture}
    \caption{1D grid enhanced with ghost cells beyond the domain boundaries, indicated by the hollow circles.}
    \label{fig:ghost-cell}
\end{figure}

The ghost-cell method enhances the discretization with ghost cells beyond the domain boundary $\partial \Omega$ (with domain $\Omega = [a,b]$), as displayed in Figure \ref{fig:ghost-cell}. Here we present the implementation for the fine grid.

The inflow (Dirichlet) BC is given by $u(x=a,t) = \alpha(t)$. Based on this, we compute ghost value $\text{u}_0$ as
\begin{equation}\label{eq:left_BC}
    u(x=a,t) = \alpha(t) = \frac{\text{u}_1+\text{u}_0}{2} \quad \rightarrow \quad \text{u}_0 = 2\alpha(t) - \text{u}_1.
\end{equation}
This corresponds to the first ghost cell outside the left boundary, see Figure \ref{fig:ghost-cell}. 
For the outflow BC we use a symmetric BC at the right boundary, given by $\frac{\partial u}{\partial x}|_{{x=b},t}=0$. This is implemented by taking $\text{u}_{N+1}=\text{u}_{N}$, where $\text{u}_{N+1}$ corresponds to the first ghost cell outside the right boundary, see Figure \ref{fig:ghost-cell}.

\subsection{Implementation for the filtered system}

\begin{figure}[ht]
    \centering
        \includegraphics[width = 0.50\textwidth]{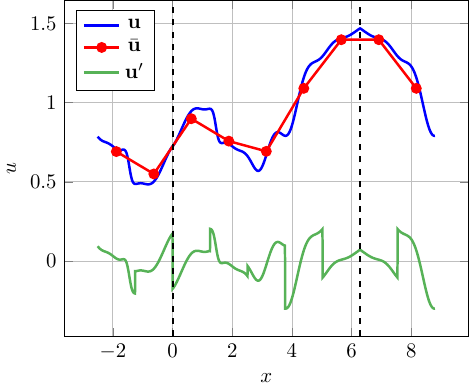}
    \caption{An example solution $\mathbf{u}$ with $N=1000$ filtered onto a coarse grid with $I=5$ extended past $\partial\Omega$, according to \eqref{eq:ref_left} ($\alpha = \frac{7}{10}$) for the left boundary and \eqref{eq:ref_right} for the right boundary. $\partial\Omega$ is indicated by the dashed vertical lines. $\bar{\mathbf{u}}$ is extended past $\partial\Omega$ according to \eqref{eq:filtered_left} and \eqref{eq:filtered_right} and $\mathbf{u}^\prime$ is extended according to \eqref{eq:SGS_left} and \eqref{eq:SGS_right}.}
    \label{fig:ex_BC}
\end{figure}

Our structure-preserving closure modeling framework is effectively a nonlinear stencil, due to the presence of a convolutional neural network. It therefore takes information from $k$ neighboring grid cells on each side. This means we require $k$ ghost cells on either side of the domain boundary $\partial \Omega$.
To find appropriate choices for the ghost values $\bar{\text{u}}_i$ and $\text{s}_i$  ($i = -k+1,\ldots,0,I+1,\ldots,I+k$) we consider the fine-grid solution $\mathbf{u}$ and appropriately extend this past the domain boundary, see Figure \ref{fig:ex_BC}. 

\subsubsection{Inflow BC}

For the left inflow BC we extend \eqref{eq:left_BC} to
\begin{equation}\label{eq:ref_left}
    \text{u}_{-i+1} = 2 \alpha(t) - \text{u}_i,\qquad i=1,2,\ldots
\end{equation}
We can rewrite this as a function of $\bar{\text{u}}_i$ and SGS content $\boldsymbol{\mu}_i$:
\begin{equation}
    (\bar{\text{u}}_{-i+1} + \mu_{-i+1,J-j}) = 2 \alpha(t) - (\bar{\text{u}}_{i} + \mu_{i,1+j}),\qquad 1\leq i\leq k,\quad 0 \leq j \leq J-1.
\end{equation}
This can be split into a filtered part:
\begin{equation}\label{eq:filtered_left}
    \bar{\text{u}}_{-i+1}  = 2\alpha(t) - \bar{\text{u}}_{i} ,\qquad 1\leq i\leq k,
\end{equation}
which yields the ghost values for $\bar{\mathbf{u}}$ past the left boundary, and a SGS part:
\begin{equation}
    \mu_{-i+1,J-j} =  - \mu_{i,1+j},\qquad 1\leq i\leq k,\quad 0 \leq j \leq J-1.
\end{equation}
The latter can be simplified as
\begin{equation}\label{eq:SGS_left}
    \boldsymbol{\mu}_{-i+1} = -\mathbf{P}\boldsymbol{\mu}_i,\qquad 1\leq i\leq k,
\end{equation}
where $\mathbf{P} \in \mathbb{R}^{J\times J}$ is the permutation matrix that represents the reflection across the boundary.

\subsubsection{Outflow BC}

For the symmetric outflow BC we extend the fine-grid solution past the domain as
\begin{equation}\label{eq:ref_right}
\text{u}_{N+i}=\text{u}_{N-i+1},\qquad i=1,2,\ldots
\end{equation}
In terms of $\bar{\text{u}}_i$ and $\boldsymbol{\mu}_i$ this becomes
\begin{equation}
    \bar{\text{u}}_{I+i} + \mu_{I+i,1+j}=\bar{\text{u}}_{I-i+1} + \mu_{I-i+1,J-j}  ,\qquad 1\leq i\leq k,\quad 0 \leq j \leq J-1,
\end{equation}
which can again be split into an equation for the ghost values for $\bar{\mathbf{u}}$:
\begin{equation}\label{eq:filtered_right}
    \bar{\text{u}}_{I+i}=\bar{\text{u}}_{I-i+1},\qquad 1\leq i\leq k,
\end{equation}
and a SGS part
\begin{equation}\label{eq:SGS_right}
\boldsymbol{\mu}_{I+i} = \mathbf{P}\boldsymbol{\mu}_{I-i+1}  ,\qquad 1\leq i\leq k.
\end{equation}

\subsection{BCs for the SGS variables}\label{sec:symm_BCs_s}

After extending the solution using the ghost-cell method we aim to find the value for $\mathbf{s}$ beyond the boundary. In order to find these ghost values we make use of the fact that we can interpret the operation $\mathbf{t}^T\mathbf{t}$ as the back and forth projection from a reduced basis to the physical SGS basis and back. This means that 
\begin{equation*}
    \boldsymbol{\mu}_i \approx \sqrt{J} \mathbf{t} \text{s}_i
\end{equation*}
and 
\begin{equation*}
    \text{s}_i = J\mathbf{t}^T\mathbf{t} \text{s}_i.
\end{equation*}
The idea is to project $\text{s}_i$ back onto physical SGS space, apply the $\mathbf{P}$ operator, and then project back to find the appropriate ghost values. Using this idea we obtain the following relations for the left and right boundary, respectively:
\begin{align}
     \boldsymbol{\mu}_{-i+1} = -\mathbf{P}\boldsymbol{\mu}_i \quad &\rightarrow \quad \text{s}_{-i+1} = -J\mathbf{t}^T\mathbf{P}\mathbf{t}\text{s}_i,\qquad 1\leq i\leq k, \\
     \boldsymbol{\mu}_{I+i} = \mathbf{P}\boldsymbol{\mu}_{I-i+1} \quad &\rightarrow \quad \text{s}_{I+i} = J\mathbf{t}^T\mathbf{P}\mathbf{t}\text{s}_{I-i+1},\qquad 1\leq i\leq k.
\end{align}
Note that $\mathbf{t}^T\mathbf{P}\mathbf{t}$ is a scalar which only needs to be computed once.

\begin{figure}[ht]
    \centering
        \includegraphics[width = 0.9\textwidth]{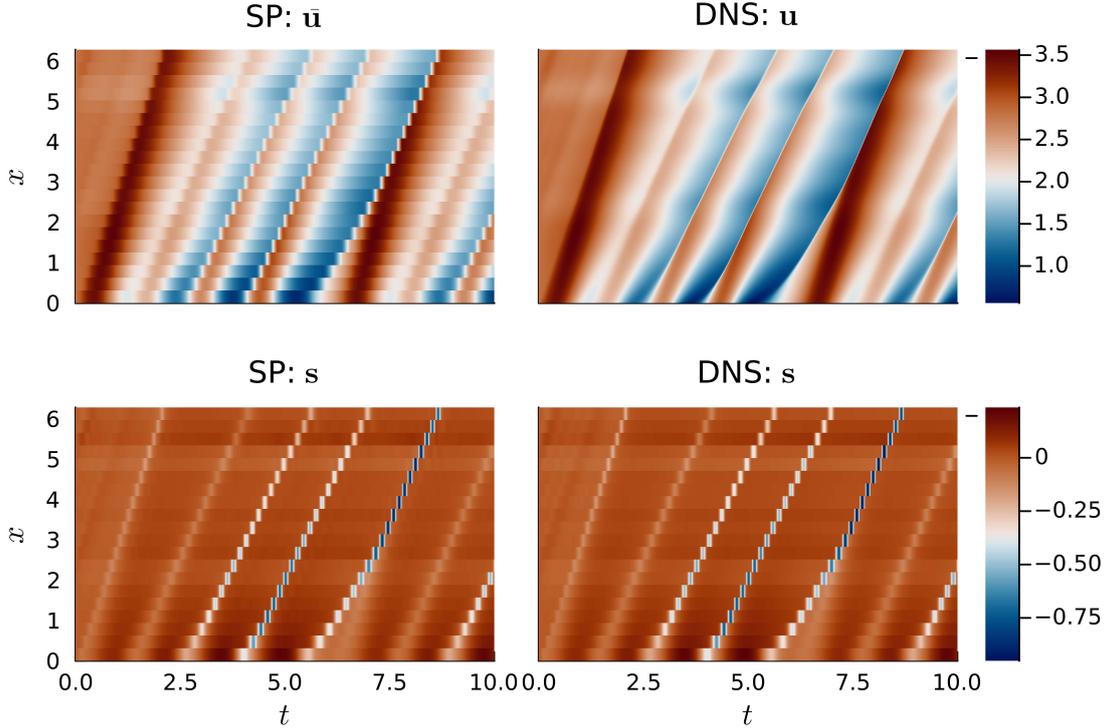}
    \caption{A simulation of Burgers' equation with I/O BCs and forcing using our trained structure-preserving closure model for $\text{DOF} = 40$ (left), along with the DNS solution for $N=1000$ (right).}
    \label{fig:sim_BC}
\end{figure}

A simulation for this inflow/outflow (I/O) BC implementation is shown in Figure \ref{fig:sim_BC}. Here we simulate Burgers' equation with I/O BCs, and some additional forcing, see \ref{sec:training}. We depict the solution produced by our structure-preserving closure modeling framework and compare it to the DNS.
\section{Training procedure}\label{sec:training}

In order to train our machine learning-based closure models we first require DNSs to generate reference data. This reference data serves as a target for our machine learning models to reproduce. In this section we describe how we randomly generate simulation conditions (initial conditions, BCs, and forcing) for closure model training and testing. In addition, we describe the training procedure and the chosen hyperparameter values, obtained from the hyperparameter tuning procedure.

\subsection{Generating training data}

To generate initial conditions, forcing, and unsteady Dirichlet BCs we make use of the following parameterized Fourier decomposition (with parameters $\alpha_1,\alpha_2,\alpha_3 \in \mathbb{R}$): 
\begin{equation}\label{eq:fourier_init}
    \xi(y;\alpha_1,\alpha_2,\alpha_3) = \alpha_1+\frac{\alpha_2}{\sqrt{M}}\sum^M_{i=2}C_{i1}\sin\left(i\frac{2\pi}{\alpha_3}y\right) + C_{i2}\cos\left(i\frac{2\pi}{\alpha_3}y\right),
\end{equation}
where $M$ is uniformly sampled from $\{2,3,\ldots,8\}$ and $C_{ij}\sim p$ from
\begin{equation}
    p(y) = \begin{cases}
        1, & \text{for } \frac{1}{2}\leq |y| \leq 1,\\
        0, & \text{elsewhere}.
        \end{cases}
\end{equation}
In the case of Burgers' equation we carry out 100 reference simulations on a uniform grid with $N=1000$ on the domain $\Omega = [0,2\pi]$ for $\nu = 0.01$. To march the solution forward in time we employ an RK4 scheme with a time step size of $\Delta t=2.5\times10^{-3}$ and simulate up to $T=10$ \cite{RK4_Butcher:2007}. 50 simulations are carried out using periodic BCs and 50 with I/O BCs. For the periodic case the initial condition is given by $u(x,t=0)=\xi(x;2,1,|\Omega|)$. For the I/O case the inflow condition is given by $u(0,t) = \xi(t,2,1,2\pi)$ and the outflow condition by a symmetric BC on the right side of the domain. The implementation of the BCs is described in \ref{sec:BCs}. The initial condition is given by a constant valued function, equal to the inflow condition at $t=0$. In addition, we also add a steady forcing term $F(x) = \xi(x;0,\frac{1}{2},|\Omega|)$ to the RHS of \eqref{eq:burgerseq} for the I/O case. 

With regards to the KdV equation we employ a uniform grid with $N=600$ on the domain $\Omega = [0,32]$ for $\varepsilon=6$ and $\mu = 1$. The solution is marched forward in time using an RK4 scheme with a time step size of $\Delta t=10^{-4}$, up to $T=10$. In this case we only consider periodic BCs, with the initial condition given by $u(x,0)=\xi(x;0,\frac{3}{5},|\Omega|)$, and perform 100 reference simulations.

For both Burgers' and KdV reference data is saved at each time interval of $5\times10^{-3}$. We randomly sample 10\% of the data from these datasets to generate the two datasets used for training (one for each equation). Both of these are split into a training (70\%) and validation set (30\%).
For testing purposes the unseen simulation conditions are generated in a similar manner, but with different randomly sampled $M$ and $C_{ij}$.

\subsection{Hyperparameters and tuning}

\begin{figure}[ht]
    \centering
    \begin{subfigure}[b]{0.466\textwidth}
        \includegraphics[width = \textwidth]{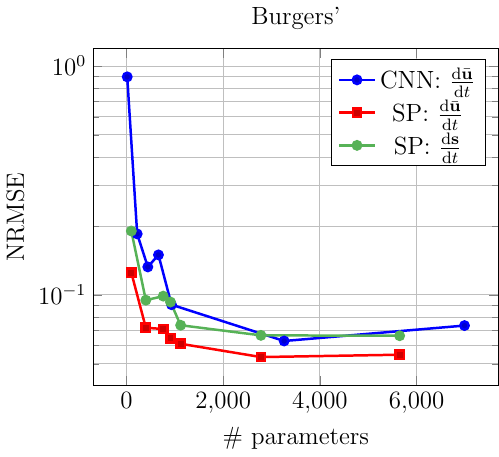}
    \end{subfigure}
    \begin{subfigure}[b]{0.48\textwidth}
        \includegraphics[width = \textwidth]{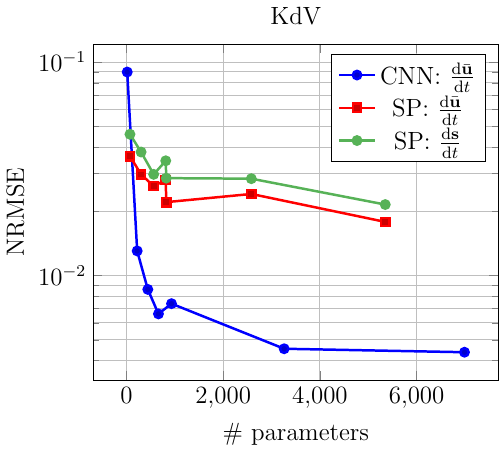}
    \end{subfigure}
    \caption{NRMSE for reproducing the RHS for each of the considered hyperparameter configurations for Burgers' (left) and KdV (right) averaged over the validation set for $\text{DOF} = 60$.}
    \label{fig:HPtuning}
\end{figure}

The chosen hyperparameters for our SP closure and the vanilla CNN are displayed in Table \ref{tab:hyperparameters}. The weights and biases are initiated using the Glorot normal initialization algorithm \cite{pmlr-v9-glorot10a}. They are optimized using the Adam optimization algorithm with parameters $\alpha$ (learning-rate), $\beta_1$ (decay rate for the first momentum estimates), $\beta_2$ (decay rate for the second momentum estimates), $\epsilon$ (small constant to combat numerical instability) \cite{kingma_adam}.
Hyperparameters are selected based on how well the trained closure models reproduce the RHS for the solution snapshots present in the validation set, corresponding to $\text{DOF}=60$. For this purpose, models are trained without trajectory fitting. For the hyperparameter optimization we opt to vary the number of hidden layers, for which we consider $\{0,1,2\}$, and the number of channels per hidden layer, for which we consider $\{10,20,30\}$.
The performance of each of the trained closure models is shown in figure \ref{fig:HPtuning}.
The best performing combination of hyperparameters (displayed in Table \ref{tab:hyperparameters}), for each equation, is selected to train the final closure models. This time trajectory fitting is included.
\begin{table}[h]
\centering 
\caption{\label{tab:hyperparameters} Hyperparameters for the trained closure models}
\begin{tabular}{|c|c|c|}
\hline
\textbf{hyperparameter}                & \textbf{CNN}           & \textbf{SP}            \\ \hline
$\alpha$                               & $10^{-3}$              & $10^{-3}$              \\ \hline
$\beta_1$                     & 0.9                    & 0.9                    \\ \hline
$\beta_2$                     & 0.999                  & 0.999                  \\ \hline
$\epsilon$                             & $10^{-8}$              & $10^{-8}$              \\ \hline
mini-batch size                    & 20                     & 20                     \\ \hline
\# iterations derivative fitting & 100                    & 100                    \\ \hline
\# iterations trajectory fitting & 20                     & 20                     \\ \hline
trajectory fitting $\overline{\Delta t}$ (Burgers')         & 0.01                   & 0.01                   \\ \hline
trajectory fitting $\overline{\Delta t}$ (KdV)         & $5\times10^{-3}$                    & $5\times10^{-3}$                   \\ \hline
trajectory fitting \# time steps (Burgers')        & 5                   & 5                   \\ \hline
trajectory fitting \# time steps (KdV)        & 20                   & 20                   \\ \hline
nonlinear activation function (underlying) CNN           & ReLU                   & ReLU                   \\ \hline
final activation function (underlying) CNN           & linear                   & linear                  \\ \hline
kernel size                   & 7                      & 5                      \\ \hline
stride                       & 1                      & 1                      \\ \hline
\# hidden layers              & 2                      & 2                      \\ \hline
\# channels per hidden layer Burgers'  & 20                     & 20                     \\ \hline
 \# channels per hidden layer KdV  & 20                     & 30                     \\ \hline
total \# parameters Burgers'         & 3261                  & 2780                   \\ \hline
total \# parameters KdV             &   3261                & 5352                  \\ \hline
$B$ (Burgers')                             & - & 1 \\ \hline
$B$ (KdV)                             & - & 2 \\ \hline
\end{tabular}
\end{table}

\end{document}